\documentclass[12pt]{amsart}
\usepackage{latexsym, amssymb, amsfonts, amscd}
\newenvironment{pf}{\proof[\proofname]}{\endproof}
\theoremstyle{plain}
\newtheorem{Thm}{Theorem}[section]
\newtheorem{Cor}[Thm]{Corollary}

\newtheorem{Prop}[Thm]{Proposition}
\newtheorem{Lem}[Thm]{Lemma}
\newtheorem{Cl}[Thm]{Claim}
\newtheorem{Conj}[Thm]{Conjecture}

\theoremstyle{definition}
\newtheorem{Def}[Thm]{Definition}

\newtheorem{Rem}[Thm]{Remark}

\newtheorem{Emp}[Thm]{}

\newtheorem{Not}[Thm]{Notation}

\numberwithin{equation}{section}

\newcommand{\pp}{\boxtimes}

\newcommand{\B}[1]{\mathbb#1}
\newcommand{\ov}{\overline}
\newcommand{\cal}[1]{\mathcal{#1}}
\newcommand{\C}[1]{\cal#1}

\newcommand{\isom}{\overset {\thicksim}{\to}}

\newcommand{\om}{\omega}

\newcommand{\hra}{\hookrightarrow}
\newcommand{\wt}{\widetilde}
\newcommand{\Gm}{\Gamma}
\newcommand{\G}{\bo{G}}
\newcommand{\bo}{\mathbf}

\newcommand{\g}{\C{G}}

\newcommand{\Dt}{\Delta}
\newcommand{\bs}{\backslash}

\newcommand{\m}{^{\times}}

\newcommand{\al}{\alpha}

\newcommand{\la}{\lambda}

\newcommand{\rl}[1]{Lemma \ref{L:#1}}
\newcommand{\rn}[1]{Notation \ref{N:#1}}
\newcommand{\rcl}[1]{Claim \ref{C:#1}}
\newcommand{\rp}[1]{Proposition \ref{P:#1}}
\newcommand{\rr}[1]{Remark \ref{R:#1}}

\newcommand{\re}[1]{\ref{E:#1}}
\newcommand{\rco}[1]{Corollary \ref{C:#1}}

\newcommand{\rt}[1] {Theorem \ref{T:#1}}
\newcommand{\rd}[1]{Definition \ref{D:#1}}
\newcommand{\sm}{\smallsetminus}

\newcommand{\ssc}{\operatorname{sc}}
\newcommand{\ab}{\operatorname{ab}}
\newcommand{\pure}{\operatorname{pure}}
\newcommand{\har}{\operatorname{char}}
\newcommand{\der}{\operatorname{der}}
\newcommand{\pr}{\operatorname{pr}}
\newcommand{\cusp}{\operatorname{cusp}}
\newcommand{\Ker}{\operatorname{Ker}}

\newcommand{\Spec}{\operatorname{Spec}}

\newcommand{\IC}{\operatorname{IC}}

\newcommand{\Frob}{\operatorname{Frob}}

\newcommand{\Lie}{\operatorname{Lie}}
\newcommand{\ad}{\operatorname{ad}}
\newcommand{\rk}{\operatorname{rk}}
\newcommand{\Id}{\operatorname{Id}}

\newcommand{\red}{\operatorname{red}}

\newcommand{\Hom}{\operatorname{Hom}}
\newcommand{\Ext}{\operatorname{Ext}}
\newcommand{\Isom}{\operatorname{Isom}}
\newcommand{\fq}{\B{F}_q}
\newcommand{\fqm}{\B{F}_{q^m}}
\newcommand{\ql}{\B{Q}_l}
\newcommand{\ga}{G(\B{A})}
\newcommand{\qlbar}{\overline{\ql}}

\begin{document}


\title[Moduli spaces of principal $F$-bundles]%
{Moduli spaces of principal $F$-bundles}
\author[Yakov Varshavsky]{Yakov Varshavsky}
\address{Institute of Mathematics, Hebrew University, Givat Ram, 91904 Jerusalem, 
ISRAEL}
\email{vyakov@math.huji.ac.il}
\thanks{This research was supported by ISF (grant No. 38/01-1)}

\begin{abstract}
In this paper we construct certain moduli spaces, which we call moduli spaces of (principal) 
$F$-bundles, and study their basic properties. These spaces are associated to triples 
consisting of a smooth projective geometrically connected curve over a finite field, a split 
reductive group $G$, and an irreducible algebraic representation $\ov{\om}$ of 
$(\widehat{G})^n/Z(\widehat{G})$. Our spaces generalize moduli spaces of $F$-sheaves, studied
by Drinfeld and Lafforgue, which correspond to the case $G=GL_r$ and $\ov{\om}$ is the tensor 
product of the standard representation and its dual. 
The importance of the moduli spaces of $F$-bundles is due to the belief that Langlands correspondence
is realized in their cohomology.
\end{abstract}
\maketitle


\section{Introduction}
Let $X$ be a smooth projective curve geometrically connected over
a finite field $\B{F}_q$, $F=\B{F}_q(X)$ the field of
rational functions on $X$, $\B{A}=\B{A}_F$ the ring of
adeles of $F$, $\Gm_F$ the absolute Galois group of
$F$, and $l$ a fixed prime, not dividing $q$.

Recall that the Langlands correspondence for $GL_r$, proved by Lafforgue 
(\cite{La}), associates  an irreducible $\ell$-adic representation
$\rho_{\pi}:\Gm_F\to GL_r(\qlbar)$ to every
cuspidal representation $\pi$ of $GL_r(\B{A})$, whose central character is
of finite order. As a result, to each pair consisting of 
$\pi$ and an algebraic representation $\om$ of
$GL_r$, Langlands correspondence associates an $\ell$-adic representation
$\rho_{\pi,\om}:=\om\circ\rho_{\pi}$ of $\Gm_F$.

Let $G$ be a split
reductive group over $\B{F}_q$, hence over $F$, and let
$\widehat{G}/\qlbar$ be the dual group of $G$.
Then Lafforgue theorem together with Langlands functoriality conjecture
predicts that to every pair $(\pi,\om)$
consisting of a tempered cuspidal representation $\pi$ of $G(\B{A})$ with 
finite order central character and an algebraic representation $\om$ of 
$\widehat{G}$, one can associate an $\ell$-adic representation $\rho_{\pi,\om}$ 
of $\Gm_F$, whose $L$-function equals 
that of $(\pi,\om)$.

More generally, for each $n\in \B{N}$ let $F^{(n)}=\B{F}_q(X^n)$ be the field of rational
functions of $X^n$. Then $\pi$ together with
an $n$-tuple $\overline{\om}=(\om_1,\ldots,\om_n)$ of representations of
$\widehat{G}$ give rise to an $\ell$-adic representation
$\rho_{\pi,\overline{\om}}$ of  $\Gm_{F^{(n)}}$, defined as the
composition of the natural restriction map
$\Gm_{F^{(n)}}\to (\Gm_F)^n$ with representation
$\boxtimes_{i=1}^n \rho_{\pi,\om_i}$ of $(\Gm_F)^{n}$. In particular, the 
Langlands conjecture associates to every pair consisting of 
$\pi$ and an irreducible representation $\overline{\om}$ of $(\widehat{G})^n$ a certain
$\ell$-adic representation $\rho_{\pi,\overline{\om}}$ of $\Gm_{F^{(n)}}$.

Furthermore, it is generally believed that the corresponding
$\rho_{\pi,\overline{\om}}$ is ``motivic'', that is, there exists
an algebraic variety $\cal X_{\pi,\overline{\om}}/F^{(n)}$ such
that $\rho_{\pi,\overline{\om}}$ is a subquotient of its cohomology. 
Moreover, for certain $\ov{\om}$'s one hopes to find 
an algebraic ``space'' $\cal X_{\ov{\om}}$ which ``realizes all of $\rho_{\pi,\overline{\om}}$''. 
By this we mean that $\cal X_{\ov{\om}}$ is equipped with an action $\ga$, and 
$H_c^*(\cal X_{\ov{\om}}, \IC(\ov{\B{Q}_l}))$ has a ``motivic''
subquotient isomorphic to the direct sum of the $(\pi\pp
\rho_{\pi,\overline{\om}})$'s, taken with certain multiplicities.

When $G=GL_r$, $n=2$ and $\ov{\om}$ is the product of the standard representation
and its dual, the existence of $\C{X}_{\ov{\om}}$ was proved by Lafforgue (\cite{La}),
generalizing an earlier work of Drinfeld (\cite{Dr1, Dr}). On the other hand, 
the required space cannot exist in the case $G=GL_2$, $n=1$ and
$\om$ is the standard representation. Indeed, the subquotient should be defined over 
$\B{Q}_l$, but the direct sum $\bigoplus_{\pi}(\pi\pp \rho_{\pi})$ is not 
(see \cite{Ka}). Thus  $\C{X}_{\ov{\om}}$ can exist only for certain $\ov{\om}$'s.

The goal of this paper is to construct a candidate of 
$\C{X}_{\ov{\om}}$ for each irreducible representation $\overline{\om}$ of 
$(\widehat{G})^n/Z(\widehat{G})$ (where $Z(\widehat{G})$ is the center of $\widehat{G}$, 
embedded diagonally in $(\widehat{G})^n$) and to study its basic properties.
By analogy with $F$-sheaves, introduced by Drinfeld, we will call our spaces
moduli spaces of (principal) $F$-bundles.

More precisely, for each $n\in\B{N}$ we construct a ``space''
$\cal X_n$ over $F^{(n)}$, equipped with an action of $\ga$. Next
for each irreducible representation $\overline{\om}$ of $(\widehat{G})^n/Z(\widehat{G})$
we construct a $\ga$-invariant closed ``subspace'' $\C{X}_{\ov{\om}}$. 
Then we construct a ``cuspidal''
subquotient of $H_c^0(\cal X_{\ov{\om}}, \IC(\ov{\B{Q}_l}))$, in
which Langlands correspondence ``should be realized''. This is especially plausible in the case
$G=GL_r$ (see Conjecture
\ref{Conj}), when Langlands correspondence is known, so  
the question is well posed. 
As evidence, we show that our conjecture holds in
the Drinfeld's case and holds ``up to $r$-negligibles'' in 
the Lafforgue case.

Roughly speaking, our construction can be described as follows: the
space $\C{X}_n$ classifies triples consisting of an
$n$-tuple $(x_1,\ldots,x_n)\in X^n$, a $G$-bundle $\C{G}$ on $X$,
and an isomorphism $\phi$ between the restrictions of $\C{G}$ and its
Frobenius  twist ${}^{\tau}\C{G}$ to the complement 
of (the graphs of) the $x_i$'s. 

To define $\C{X}_{\ov{\om}}$'s, observe that each irreducible representation $\ov{\om}$ of 
$(\widehat{G})^n$ (hence of $(\widehat{G})^n/Z(\widehat{G})$) corresponds to a certain 
$n$-tuple  $(\om_1,\ldots,\om_n)$ of dominant coweights of $G$ (see \rr{corr}).
Then we define $\C{X}_{\ov{\om}}$ to be the closed substack of $\C{X}_n$, consisting of 
those triples $(\C{G};x_1,\ldots,x_n;\phi)$, for which the relative position of 
$\phi(\C{G})$ and $^{\tau}\C{G}$ at $x_i$ is less than or equal to $\om_i$ for each $i$.

By construction, $\C{X}_n$ is just a twisted version of the global affine grassmannian over 
$X^n$. Moreover, if we denote by $\C{F}_{\overline{\om}}/\C{X}_n$ the extension by zero of the
IC-sheaf of $\C{X}_{\overline{\om}}$, then the correspondence
$\overline{\om}\mapsto\C{F}_{\overline{\om}}$ is just a twisted
version of the geometric Satake correspondence (see \rt{locisom}
and \rco{strata}).

Finally, to get the required subquotient of $H^0_c(\C{X}_{\ov{\om}},\IC(\ov{\B{Q}_l}))$ 
we proceed in two steps: first, we consider its maximal pure quotient
of weight zero  $H^0_{c,\pure}(\C{X}_{\ov{\om}},\IC(\ov{\B{Q}_l}))$, and then take the subspace of 
 $H^0_{c,\pure}(\C{X}_{\ov{\om}},\IC(\ov{\B{Q}_l}))$, consisting of all elements, vanishing on the 
locus of reducible $F$-bundles (i.e., those $F$-bundles which has a $\phi$-invariant parabolic structure). 

We would like to note that our construction is a rather straightforward combination of the 
original Drinfeld construction of the moduli of $F$-sheaves, Beilinson--Drinfeld construction of 
the Hecke stacks, and geometric Satake correspondence. 
In particular, it was known to Drinfeld and some others.

\vskip 6truept
\centerline{\bf Notation and conventions}
\vskip 6truept

1) Let $G$ be a split reductive group over a finite
field $\fq$, let $G^{\der}$ be the derived group of $G$, $G^{\ssc}$ the simply-connected 
cover of $G^{\der}$, $G^{\ab}:=G/G^{\der}$ the abelinization of $G$, and $G^{\ad}$ the 
adjoint group of $G$. Let $B\supset T\subset Z$ be a Borel subgroup, a maximal torus,
and the center of $G$, respectively. We denote by $B^{\ssc}\supset T^{\ssc}\subset Z^{\ssc}$ 
the corresponding objects of $G^{\ssc}$, and similarly for $G^{\der}$ and $G^{\ad}$.

2) Let $\rho$ be the half-sum of all positive coroots of $G$.

3) By a {\em quasi-fundamental weight} of $G$ we mean the
smallest positive multiple of a fundamental weight of $G^{\ssc}$,
which belongs to $X^*(T^{\ad})\subset X^*(T)$.

4) Let $X^*_+(T)$ and $X^+_*(T)$ be the sets of dominant weights
and coweights of $G$, respectively.

5) Weights (resp. coweights) of $G$ we equip with (standard) ordering: 
$\la_1\leq\la_2$ if and only if the difference $\la_2-\la_1$ is a positive 
integral linear combination of simple roots (resp. coroots) of $G$.

6) For an algebraic group $H$, an $H$-bundle $\C{H}$ on $Y$, and a representation $V$ of
$G$, we denote the vector bundle $H\bs[\C{H}\times V]$ on $Y$ by
$\C{H}_V$. 

7) For a dominant weight $\la$ of $G$, we denote by $V_{\la}$ the 
Weyl module of $G$ with the highest weight $\la$. Also for a $G$-bundle $\C{G}$ on $Y$, 
we denote $\C{G}_{V_{\la}}$ by $\C{G}_{\la}$.

8) For a finite scheme $D$ over a field $k$, put
$\C{O}_D:=k[D]$ and $|D|:=\dim_k\C{O}_D$.

9) For a finite scheme $D$ over a field $k$, denote the Weil restriction of scalars
$R_{D/k}G$ by $G_{D}$. In particular, $G_D(k)=G(\C{O}_D)$.
More generally, for every closed embedding between finite schemes $D_1\hra D_2$, we denote  
by $G_{D_1,D_2}$ the kernel of the natural homomorphism $G_{D_2}\to G_{D_1}$.

10) For a closed point $v$ of a curve $X$, let $\C{O}_v$ and $F_v$
be the completions at $v$ of the stalk at $v$ of the structure
sheaf and the field of its fractions, respectively.

11) For an $S$-point $x$ of a scheme $X$, let $\Gm_x\subset X\times
S$ be the graph of $x$.

12) By $\Dt\subset X^n$ we denote the set of all $n$-tuples
$(x_1,\ldots,x_n)$ for which there exist $i\neq j$ with $x_i=x_j$.

13) By an IC-sheaf on a stack $Y$, we will mean the intermediate extension 
of the constant perverse $\overline{\B{Q}_l}$-sheaf on a open dense 
substack $Y^0$ of $Y$ such that the corresponding reduced stack $(Y^0)_{\red}$ is smooth.
The IC-sheaf is normalized so that it is pure of weight zero. 
The IC-sheaf on $Y$  will denote by $\IC_Y$ or simply by $\IC$.

14) For a stack $Y$ over a finite field $\fq$, we denote by $\Frob_q:Y\to Y$ the absolute Frobenius 
morphism over $\fq$.

\vskip 6truept
\centerline{\bf Acknowledgments}
\vskip 6truept

First, the author thanks V. Drinfeld, who
suggested to study moduli spaces of $F$-sheaves and mentioned that
this construction can be extended to arbitrary groups. Secondly,
the present work would not be possible without numerous
conversations with D. Gaitsgory, who explained many things to me about
$B$-structures, affine grassmannians and geometric Satake correspondence.
He also read preliminary drafts of the paper, and his suggestions significantly 
improved the exposition.  Also I thank J. Arthur, J. Bernstein and D. Kazhdan for  
their interest, stimulating conversations and useful remarks.

Different parts of the work were done while the author visited the
University of Toronto, IHP and IHES, which I thank for stimulating atmosphere and 
financial support. 
\section{Main constructions and results} \label{S:definitions}

\begin{Not} \label{N:basic}

a) Let $X$ be a smooth projective curve geometrically connected
over a field $k$, and let $Bun=Bun_G$ be the stack classifying
$G$-bundles on $X$, i.e., $Bun_G(S)=\{G-\text{bundles on }X\times
S\}$ for each scheme $S$ over $k$. 

More generally, for each finite subscheme $D\subset X$, let $Bun_D=Bun_{G,D}$ be
the stack over $Bun$ classifying $G$-bundles on $X$  with
$D$-level structures, i.e.,
$$Bun_{G,D}(S)=\{\cal G\in Bun_G(S),\psi:\cal G_{|D\times S}\isom G\times D\times S\}.$$

b) For each $\mu\in X_*(T^{\ssc})\otimes \B{Q}=X_*(T^{\ad})\otimes
\B{Q}$, let $Bun_G^{\leq \mu}$ be a substack of $Bun_G$,
consisting of $G$-bundles, whose degree of instability is bounded
by $\mu$, i.e.,
\[
Bun_G^{\leq \mu}(S)=\{\g\in Bun_G(S)|\text { for each geometric point }
{s}\in S,\text{ each }
\]
\[
 B-\text{structure } \C{B} \text{ of } \g_{{s}} \\ \text { and each }
\la\in X^*_+(T^{\ad}): \deg \C{B}_{\la}\leq \langle \mu,\la\rangle \},
\]
\noindent where $\C{B}_{\la}$ is the corresponding line bundle. 
A substack $Bun_G^{\leq \mu}$ is open in $Bun_G$ (see
\rl{open}). More generally, for every $D$ we will denote by  $Bun_{G,D}^{\leq \mu}$ the preimage of 
 $Bun_G^{\leq \mu}$ in $Bun_{G,D}$.
\end{Not}

The proof of the following basic fact will be recalled in  \re{conn}.

\begin{Lem} \label{L:conn}
The set of connected components $\pi_0(Bun_G)$ of $Bun_G$ is
canonically isomorphic to $\pi_1(G):=X_*(T)/X_*(T^{\ssc})$ (which in its turn is
canonically isomorphic to the group of characters of $Z(\widehat{G})$).
\end{Lem}

\begin{Not}
Denote by $\pi_0$ the canonical map $Bun_G\to\pi_0(Bun_G)=\pi_1(G)$, and
denote by $[\om]\in\pi_1(G)$ the class of $\om\in X_*(T)$.
\end{Not}

Next will introduce  affine grassmannians and Hecke stacks 
following Beilinson and Drinfeld
(\cite{BD}).

\begin{Def} \label{D:hecke}
a) For each $n\in\B{N}$ and each finite (not necessary non-empty) subscheme $D\subset X$, let
$Hecke_{D,n}$ be the stack, which for each scheme
$S$ over $k$, classifies triples:

i) $(\C{G},\psi),(\C{G}',\psi')\in Bun_{D}(S)$;

ii) $n$ points $x_1,x_2,\ldots,x_n\in (X\sm D)(S)$;

iii) isomorphism $\phi:\C{G}_{|(X\times S)\sm(\Gm_{x_1}\cup\ldots\cup
\Gm_{x_n})}\isom \C{G}'_{|(X\times S)\sm (\Gm_{x_1}\cup\ldots\cup
\Gm_{x_n})}$, preserving $D$-level structures (that is satisfying 
$\psi'\circ\phi_{|D\times S}=\psi$).

As in the case of $Bun_D$, we will omit $D$ from the notation when $D=\emptyset$.

b) For each $n$-tuple of dominant coweights $\overline{\om}=(\om_1,\ldots,\om_n)$ of
$G$, let $Hecke_{D,n,\overline{\om}}$ be the closed substack of
$Hecke_{D,n}$ defined by the condition that ``the relative position of
$\phi(\C{G})$ and $\C{G}'$ at $x_i$ is less than or equal to $\om_i$
for each $i$'' in the following sense:

iii)$_{\overline{\om}}$ $\phi(\C{G}_{\la})\subset\C{G}'_{\la}(\sum_{i=1}^n
\langle \la,\om_i\rangle \Gm_{x_i})$ for each dominant weight $\la$ of $G$;

iii$')_{\overline{\om}}$ $\pi_0(\C{G}_{{s}})-\pi_0(\C{G}'_{{s}})=[\sum_i\om_i]$ 
for each geometric point ${s}\in S$.
\end{Def}

By \rl{repr}, $Hecke_{D,n,\overline{\om}}$ is an
algebraic stack locally of finite type over $k$.

\begin{Rem} \label{R:inv}

a) Condition iii)$_{\overline{\om}}$ implies that for each
character $\la\in X^*(G)$ we get
$\phi(\C{G}_{\la})=\C{G}'_{\la}(\sum_{i=1}^n
\langle \la,\om_i\rangle \Gm_{x_i})$. Indeed, apply iii)$_{\overline{\om}}$ to
both $\la$ and $\la^{-1}$.

b) If $G^{\der}$ is simply connected, then iii$')_{\overline{\om}}$
is a consequence of iii)$_{\overline{\om}}$. Indeed, in this case,
$\pi_1(G)=X_*(G^{\ab})$. Hence condition  iii$')_{\overline{\om}}$
is equivalent to the equality $\deg(\C{G}_{\la})=\deg(\C{G}'_{\la})+\sum_{i=1}^n\langle \la,\om_i\rangle$ for every $\la\in X^*(G)=X^*(G^{\ab})$
(compare the proof of \rl{conn} in \re{conn}).
Thus the statement follows from a).

c) Stack $Hecke_{n}$ has a natural involution which interchanges
$\C{G},\phi$ with $\C{G}',\phi^{-1}$. This
involution sends $Hecke_{n,\overline{\om}}$ into
$Hecke_{n,-w_0(\overline{\om})}$, where $w_0$ is the longest
element of the Weyl group of $G$, acting on $X^*(T)^n$ diagonally.
To see this, note that $V_{-w_0(\la)}$ is
the dual of the Weyl module $V_{\la}$ and that the projection $X^*(T)\to
\pi_1(G)$ is constant on the orbits of the Weyl group.

d) One may consider a variant of the definition of  $Hecke_{n,\overline{\om}}$,
in which iii)$_{\overline{\om}}$ is replaced by an a priori stronger condition: 
$\phi(\C{G}_V)\subset\C{G}'_V(\sum_{i=1}^n\langle \xi,\om_i\rangle \Gm_{x_i})$ 
for each weight $\xi$ of $G$ and each representation $V$ of $G$ all of whose weights are 
less than or equal to $\xi$. Though the stack defined by this condition might be  
smaller than the original one, it follows from the Cartan decomposition that the corresponding 
reduced stacks coincide. However we do not know whether the same 
is true for stacks themselves (compare \rr{reduced}).
\end{Rem}
 
Following \cite{BD}, we also consider iterated Hecke stacks. 

\begin{Def} \label{D:ithecke}
a) For each $n\in\B{N}$ and each finite subscheme $D\subset X$, let 
$Hecke'_{D,n}$ be the stack which for each scheme
$S$ over $k$ classifies triples consisting of the following:

i) $n+1$ elements $(\C{G},\psi)=(\C{G}_0,\psi_0),(\C{G}_1,\psi_1),\ldots, (\C{G}_n,\psi_n)=(\C{G}',\psi')$
of $Bun_{D}(S)$;

ii) $n$ points $x_1,x_2,\ldots,x_n\in (X\sm D)(S)$;

iii) isomorphisms $\phi_{i}:{\C{G}_{i-1}}_{|(X\times S)\sm\Gm_{x_i}}\isom
{\C{G}_{i}}_{|(X\times S)\sm \Gm_{x_i}}$, preserving $D$-level structures, 
for each $i=1,\ldots,n$.

b) For each $n$-tuple of dominant coweights $\overline{\om}=(\om_1,\ldots,\om_n)$ of 
$G$, let $Hecke'_{D,n,\overline{\om}}$ be the closed substack of 
$Hecke'_{D,n}$ defined by the condition that each 
$$
[(\C{G}_{i-1},\psi_{i-1}),(\C{G}_i,\psi_i);x_i;\phi_i]\in Hecke_{D,1}(S)
$$
belongs to $Hecke_{D,1,\om_i}(S)$. 
\end{Def}

\begin{Rem} \label{R:ithecke}

a) Alternatively,  $Hecke'_{D,n,\overline{\om}}$ can be defined as a fiber product 
$$
Hecke_{D,1,\om_1}\times_{Bun_D} Hecke_{D,1,\om_2}\times_{Bun_D}\ldots\times_{Bun_D} Hecke_{D,1,\om_n}.
$$ 

b) We have a natural forgetful map $\pi:Hecke'_{D,n,\overline{\om}}\to Hecke_{D,n,\overline{\om}}$, 
which forgets
$(\C{G}_1,\psi_1),\ldots, (\C{G}_{n-1},\psi_{n-1})$ and replaces the 
$\phi_i$'s by their composition.
Furthermore, $\pi$ is projective (see \rl{repr}) surjective and small (see \rl{small}).

c) More generally, for each partition $n=k_1+\ldots+k_l$ one can similarly consider a partially 
iterated Hecke stack
$$
Hecke_{D,k_1}\times_{Bun_D} Hecke_{D,k_2}\times_{Bun_D}\ldots\times_{Bun_D} Hecke_{D,k_l}.
$$ 
The previously defined stacks $Hecke_{D,n}$ and  $Hecke'_{D,n}$ correspond 
to the trivial partition $n=n$ and the maximal partition $n=1+\ldots +1$, respectively.  
\end{Rem}

\begin{Not} \label{N:proj}
Denote by $p$ and $p'$ forgetful morphisms  $Hecke_{D,n}\to Bun_D$ sending the triple to 
$(\C{G},\psi)$ and $(\C{G}',\psi')$ respectively, and define by 
$Hecke^{\leq\mu}_{D,n}$ and $Hecke^{\leq\mu}_{D,n,\overline{\om}}$ the preimages of $Bun^{\leq\mu}_D$ 
under $p$. Similarly we define $Hecke'^{\leq\mu}_{D,n}$ and $Hecke'^{\leq\mu}_{D,n,\overline{\om}}$.
\end{Not}

\begin{Rem} \label{R:level}
The space $Hecke_{D,n}$ is canonically isomorphic to the restriction to $(X\sm D)^n$ of the 
fiber product $Hecke_n\times_{Bun} Bun_D$, where the map $Hecke_n\to Bun$ is either $p$ or $p'$. 
\end{Rem}

\begin{Def} \label{D:glaff}
Let $Gr_n$ (resp. $Gr_{n,\overline{\om}}$, $Gr'_n$, $Gr'_{n,\overline{\om}}$) be the stacks 
classifying the same data as $Hecke_n$ (resp
$Hecke_{n,\overline{\om}}$, $Hecke'_n$, $Hecke'_{n,\overline{\om}}$), together with a trivialization of 
$\C{G}'$. These spaces are called {\em global affine grassmannians (over $X^n$)}.
\end{Def}

Now we are ready to introduce our main object. From now on $k$ will be a
finite field $\fq$.

\begin{Not} \label{N:twist}
For a scheme $S/\fq$ and an $S$-point $A$ of a stack $\C{X}$ over $\fq$, 
we denote the $S$-point $Frob_q^*(A)$ by ${}^{\tau}A$. In particular, for a 
coherent sheaf or a $G$-bundle
$\C{F}$ over $X\times S$, we will write ${}^{\tau}\C{F}$ instead of 
$(Id_X\times Frob_{q})^*(\C{F})$.
\end{Not}

\begin{Def} \label{D:Fbundles}
For each $n\in\B{N}$ and finite subscheme $D\subset X$,
let $FBun_{D,n}$ (resp. $FBun_{D,n,\overline{\om}}$, $FBun'_{D,n,\overline{\om}}$,
$FBun^{\leq\mu}_{D,n,\overline{\om}}$) be the stack
classifying the same data i)--iii) as $Hecke_{D,n}$ (resp.
$Hecke_{D,n,\overline{\om}}$, $Hecke'_{D,n,\overline{\om}}$, 
$Hecke^{\leq\mu}_{D,n,\overline{\om}}$) together with an
isomorphism $\C{G}'\isom{}^{\tau}\C{G}$, preserving $D$-level structures.

We will call these spaces {\em moduli spaces of (principal) $F$-bundles}.
\end{Def}

\begin{Rem} \label{R:fbundles}
a) Explicitly, a stack $FBun_{D,n}$ classify triples consisting of a pair
$(\C{G},\psi)\in Bun_D(S)$,  an $n$-tuple $(x_1,\ldots,x_n)\in (X\sm D)(S)$, 
and an isomorphism 
$$
\phi:\C{G}_{|(X\times S)\sm(\Gm_{x_1}\cup\ldots\cup\Gm_{x_n})}
\isom{}^{\tau}\C{G}_{|(X\times S)\sm (\Gm_{x_1}\cup\ldots\cup\Gm_{x_n})}
$$ 
such that  
$^{\tau}\psi\circ\phi_{|D\times S}=\psi$. In particular,  $FBun_{D,n}$ is equipped with a 
natural action of the group $G(\C{O}_D)$, which replaces $\psi$ by $g\circ\psi$ 
for each $g\in G(\C{O}_D)$ and does not change all the other data.

b) In the case $G=GL_r$, $\om_1=(0,\ldots,0,-1)$ and $\om_2=(1,0,\ldots,0)$, the space of 
$F$-bundles $FBun'_{D,2,(\om_1,\om_2)}$ (resp.  $FBun'_{D,2,(\om_2,\om_1)}$) is 
the moduli space of right (resp. left) $F$-sheaves $FSh_{D,r}$ (resp. $_{D,r}FSh$)
studied by Drinfeld and Lafforgue.
\end{Rem}

\begin{Def} \label{D:adm}
An $n$-tuple $\overline{\om}=(\om_1,\ldots,\om_n)\in X_*(T)^n$ is called {\em admissible}
if the sum of the $\om_i$'s belongs to $X_*(T^{\ssc})$.
\end{Def}

\begin{Rem} \label{R:adm}
As $Frob_q$ acts trivially on $\pi_0(Bun_G)$, we get from \rl{conn} that condition  
iii$')_{\overline{\om}}$ of \rd{hecke} implies that $\overline{\om}$ is admissible 
if $FBun_{n,\ov{\om}}$ is non-empty. Conversely, if $\ov{\om}$ is admissible,
then  condition iii$')_{\overline{\om}}$ in the definition of $FBun_{n,\ov{\om}}$
holds automatically.
\end{Rem}

The following proposition, whose proof will be given in
\re{fbundles}, summarizes basic properties of the moduli spaces of $F$-bundles, generalizing
\cite[Prop. 2.3 and 3.2]{Dr1}.

\begin{Prop} \label{P:fbundles}

a) $FBun_{D,n,\overline{\om}}$ is a Deligne--Mumford stack over $(X\sm D)^n$,
locally of finite type. Moreover, connected components of
$FBun^{\leq\mu}_{D,n,\overline{\om}}$ are quotients of
quasi-projective schemes over $(X\sm D)^n$ by finite groups.
Furthermore, these components are quasi-projective schemes, if
$|D|$ is sufficiently large relative to $\mu$.

b) Every $FBun_{D,n,\overline{\om}}$ is a finite (\'etale) Galois cover of 
$FBun_{n,\overline{\om}}\times_{X^n}(X\sm D)^n$ with Galois group
$G_D(\fq)$. In particular, for each 
$D_1\subset D_2$, $FBun_{D_2,n,\overline{\om}}$ is a finite  
Galois cover of $FBun_{D_1,n,\overline{\om}}\times_{(X\sm
D_1)^n}(X\sm D_2)^n$ with Galois group $G_{D_2,D_1}(\fq)$.

c) If $\om_1=\ldots=\om_n=0$, then  $FBun_{D,n,\overline{\om}}$ is
canonically isomorphic to the product of $(X\sm D)^n$ with a discrete
stack $Bun_D(\fq)$.

d) $FBun_{D,n,\overline{\om}}$ is non-empty if and only if $\overline{\om}$ is admissible.

e) The forgetful morphism  $\pi:FBun'_{D,n,\overline{\om}}\to FBun_{D,n,\overline{\om}}$
is projective. Moreover, $\pi$ is an isomorphism over $X^n\sm\Dt$. 
\end{Prop}

\begin{Rem} \label{R:corr}
$n$-tuples of dominant coweights of $G$ are in canonical bijection with dominant weights of 
$(\widehat{G})^n$ and hence with irreducible representations of $(\widehat{G})^n$.
Under this bijection, admissible $n$-tuples 
correspond to representations, trivial on $Z(\widehat{G})$.
\end{Rem}

\begin{Not}
Denote by $[n]$ the set $\{1,\ldots,n\}$, and for each set
$A$ we will identify $A^n$ with the set of functions on $[n]$ with values in $A$. In
particular, every map $\eta:[n]\to [k]$ induces a map
$\eta^*:A^k\to A^n$. Moreover, if $A$ is an abelian group, then 
$\eta$ induces also a map $\eta_!:A^n\to A^k$ 
(defined as $\eta_!(f)(i):=\sum_{j\in\eta^{-1}(i)}f(j)$).
\end{Not}

\begin{Emp} \label{E:str} {\bf Stratification.}
All of the above stacks have natural stratifications:

a) The stratification of $Hecke'_n$ (and hence of $FBun'_n$ and $Gr'_n$) is indexed by 
$n$-tuples $\overline{\om}\in X^+_*(T)^n$. To see this, define a partial order on 
$X^+_*(T)^n$ by the 
rule $\overline{\om}'\leq\overline{\om}$ if and only if $\overline{\om}'(i)\leq\overline{\om}(i)$ 
for each $i\in[n]$. Then $Hecke'_{n,\overline{\om}'}$ is contained in
$Hecke'_{n,\overline{\om}}$ if $\overline{\om}'\leq\overline{\om}$.
For each  $\overline{\om}\in X^+_*(T)^n$, let $Hecke'^0_{n,\overline{\om}}$ be the complement in 
$Hecke'_{n,\overline{\om}}$ of the union of all $Hecke_{n,\overline{\om}'}$'s with
$\overline{\om}'<\overline{\om}$. Then $\{Hecke'^0_{n,\overline{\om}}\}_{\ov{\om}}$ gives us a
required stratification of $Hecke'_n$, and we denote by $\{FBun'^0_{n,\overline{\om}}\}_{\ov{\om}}$ and 
$\{Gr'^0_{n,\overline{\om}}\}_{\ov{\om}}$ the induced stratifications of $FBun'_n$ and $Gr'_n$, 
respectively. 

b) The stratification of $Hecke_n$ (and hence of $FBun_n$ and $Gr_n$) is indexed by equivalence classes of triples
$(k,\eta,\overline{\om})$, where $k\leq n$ is a positive integer,
$\eta:[n]\to [k]$ is a surjection, $\overline{\om}$ is an element of $X^+_*(T)^k$, and the 
equivalence relation is given by the rule
$(k,\sigma\circ\eta,\overline{\om})\sim(k,\eta,\sigma^*(\overline{\om}))$ 
for each $\sigma\in S_k$.
The set of equivalence classes has a natural partial order defined by the rule that
$[(k',\eta',\overline{\om}')]\leq[(k'',\eta'',\overline{\om}'')]$ if
and only if there is a surjection $\eta:[k'']\to[k']$ such that
$\eta'=\eta\circ\eta''$ and $\overline{\om}'\leq\eta_!(\overline{\om}'')$.

For each $\C{T}=[(k,\eta,\overline{\om})]$, let
$Hecke_{n,\C{T}}\subset Hecke_n$ be the image
of $Hecke_{k,\eta_!(\overline{\om})}$ under the closed
embedding $Hecke_k\hra Hecke_n$, induced by $\eta^*:X^k\hra X^n$. 
Then $Hecke_{n,\C{T}'}$ is contained in
$Hecke_{n,\C{T}}$ if $\C{T}'\leq \C{T}$.
Denote by $Hecke^0_{n,\C{T}}$ the complement in 
$Hecke_{n,\C{T}}$ of the union of all
$Hecke_{n,\C{T}'}$'s with $\C{T}'<\C{T}$.  Then $\{Hecke^0_{n,\C{T}}\}_{\C{T}}$ gives us a
required stratification of $Hecke'_n$, and we denote by  $\{FBun^0_{n,\C{T}}\}_{\C{T}}$  and  
$\{Gr^0_{n,\C{T}}\}_{\C{T}}$ the induced stratifications of $FBun_n$ and $Gr_n$, 
respectively. 

For each $\ov{\om}\in X_*^+(T)$, we denote $\C{T}_{\ov{\om}}:=[(n,\Id,\ov{\om})]$ simply by $\ov{\om}$.
This will not lead to confusion, since each $Hecke_{n,\C{T}_{\ov{\om}}}$ coincides with  
$Hecke_{n,\ov{\om}}$ and since we have $\C{T}_{\ov{\om}'}\leq \C{T}_{\ov{\om}}$ if and only if 
$\ov{\om}'\leq\ov{\om}$.
\end{Emp}

The following theorem and its corollary, which will be proved in
\re{locmod} and \re{strata}, respectively, imply that locally in
the \'etale topology, $FBun_{n,\overline{\om}}$ (resp. $FBun'_{n,\overline{\om}}$) is isomorphic to 
$Gr_{n,\overline{\om}}$ (resp. $Gr'_{n,\overline{\om}}$). This result generalizes the
corresponding result of Drinfeld (\cite[Prop. 3.3]{Dr1}), asserting that
the moduli space of $F$-sheaves is smooth.

\begin{Thm} \label{T:locisom}

$Gr_{n,\overline{\om}}$ is a local model of
$FBun_{D,n,\overline{\om}}$. In other words, for every point $y\in
FBun_{D,n,\overline{\om}}$, there exists an \'etale neighborhood
$p_1:U_y\to FBun_{D,n,\overline{\om}}$ of $y$ and an \'etale
morphism $p_2:U_y\to Gr_{n,\overline{\om}}$. Moreover, $p_1$ and
$p_2$ induce the same stratification of $U_y$ and the same
morphism $U_y\to X^n$. Furthermore, $p_2$ lifts to an \'etale morphism 
$U_y\times_{FBun_{D,n,\overline{\om}}} FBun'_{D,n,\overline{\om}}\to Gr'_{n,\overline{\om}}$,
compatible with stratifications. In particular, $Gr'_{n,\overline{\om}}$ is a local model
of $FBun'_{D,n,\overline{\om}}$.
\end{Thm}

\begin{Cor} \label{C:strata}

a) The open stratum $FBun^0_{D,n,\overline{\om}}$ of $FBun_{D,n,\overline{\om}}$ 
(resp. $FBun'^0_{D,n,\overline{\om}}$ of $FBun'_{D,n,\overline{\om}}$)
is dense. It is non-empty if and only if $\overline{\om}$ is
admissible. 

b) The reduced stacks $(FBun^0_{D,n,\overline{\om}})_{\red}$ and 
$(FBun'^0_{D,n,\overline{\om}})_{\red}$ are smooth over \break
$(X\sm D)^n$ of relative dimension $\sum_{i=1}^n\langle 2\rho,\om_i\rangle$.
Furthermore, both $FBun^0_{D,n,\overline{\om}}$ and  $FBun'^0_{D,n,\overline{\om}}$ are reduced, unless
$\har  \fq=2$, and $G$ has a direct factor isomorphic to $PGL_2$ or $PO_{2m+1}$.

c) The IC-sheaf of $FBun_{D,n,\overline{\om}}$ (resp. $FBun'_{D,n,\overline{\om}}$) is the
restriction (up to a homological shift and Tate twist) of that of
$Hecke_{D,n,\overline{\om}}$ (resp. $Hecke'_{D,n,\overline{\om}}$). In particular, its restriction to
each stratum is a direct sum of complexes of the form
$\overline{\B{Q}_l}(k+n/2)[2k+n]$ with $k\in\B{Z}$. 

d) The forgetful morphism $\pi:FBun'_{D,n,\overline{\om}}\to FBun_{D,n,\overline{\om}}$
is projective, surjective and small. In particular, the intersection cohomology (with compact support) of 
$FBun_{D,n,\overline{\om}}$ coincides with that of $FBun'_{D,n,\overline{\om}}$.
\end{Cor}

\begin{Rem} \label{R:reduced1}
The second assertion of b) says that in most  cases the open strata  $FBun^0_{D,n,\overline{\om}}$ 
and $FBun'^0_{D,n,\overline{\om}}$ are reduced. However, we do not know whether the same is true
for the full stacks $FBun_{D,n,\overline{\om}}$ and $FBun'_{D,n,\overline{\om}}$ (compare \rr{reduced}).
\end{Rem}

\begin{Def} \label{D:red}
We will call an $F$-bundle $(\C{G}; x_1,\ldots,x_n;\phi)$ {\em
reducible} if there exists a maximal parabolic subgroup
$P\subset G$ and a $P$-structure $\C{P}$ of $\C{G}$ such that
$\phi$ induces a rational isomorphism between $\C{P}$ and
${}^{\tau}\C{P}$.
\end{Def}

The following result, proved in \re{red}, shows that ``at infinity''
all $F$-bundles are reducible.

\begin{Not} \label{N:max}
Let $d(\overline{\om})$ be the maximum of the $\langle \sum_{k=1}^n\om_k+4g\rho,\lambda_i\rangle $'s taken over
the set of all fundamental weights $\lambda_i$ of $G^{\ssc}$, where $g$ is the genus of $X$.
\end{Not}

\begin{Thm} \label{T:red}
Every $F$-bundle from $FBun_{n,\ov{\om}}\sm FBun_{n,\ov{\om}}^{\leq d(\overline{\om})\rho}$ is reducible.

\end{Thm}

\begin{Rem} \label{R:Behrend}
Using methods and results of K. Behrend (\cite{Be}), one can show that \rt{red} still 
remains true if $d(\ov{\om})$ is replaced by the maximum of the 
$\langle \sum_{k=1}^n\om_k,\lambda_i\rangle $'s. 
In particular, the bound $d(\ov{\om})$ can be made independent of the curve $X$. 
However, \rt{red} seems to be sufficient for all the applications, and the proof of a better bound 
is much more involved.
\end{Rem}

\begin{Not}
Let $FBun_{*,n}$ be the generic fiber over $X^n$ of the inverse limit of the 
$FBun_{D,n}$'s, and let $FBun_{*,n,\overline{\om}}$ be the corresponding closed substack.
\end{Not}


\begin{Emp} \label{E:reduc}
For each maximal parabolic $P$, let $FBun_{P,n}$ be the
stack classifying the data consisting of an $F$-bundle
$(\C{G}; x_1,\ldots,x_n;\phi)$ and a $P$-structure $\C{P}$ of
$\C{G}$ such that $\phi$ induces a rational isomorphism between
$\C{P}$ and ${}^{\tau}\C{P}$. We have a natural forgetful map
$FBun_{P,n}\to FBun_n$, whose image is the set of all reducible
$F$-bundles, corresponding to $P$. More generally, define
$FBun_{P,D,n}$, $FBun_{P,n,\overline{\om}}$ be the fiber product 
of $FBun_{P,n}$ over $FBun_n$ with $FBun_{D,n}$,
$FBun_{n,\overline{\om}}$, respectively.
\end{Emp}

\begin{Def} \label{D:orisph}
By an {\em orispheric substack} we will call the image in $FBun_{D,n,\overline{\om}}$ 
(resp.  $FBun_{*,n,\overline{\om}}$) of an irreducible component of $FBun_{P,D,n,\overline{\om}}$  
(resp.  $FBun_{P,*,n,\overline{\om}}$).
\end{Def}

A more precise version (\rp{par1}) of the following result
generalizes the corresponding results of Drinfeld (\cite[Prop.~4.3]{Dr1})
and Lafforgue (\cite[II, Thm.~5]{La1}).

\begin{Prop} \label{P:par}
Every orispheric substack of $FBun_{*,n,\overline{\om}}$ is closed.
\end{Prop}

The following simple result, proven in \re{action}, provides us
with a space over $F^{(n)}$, equipped with an action of $G(\B{A})$.

\begin{Prop} \label{P:action}

a) The group  $Z(\B{A})/Z(F)$ acts naturally on $FBun_{D,n}$ and preserves each 
$FBun^{\leq\mu}_{D,n,\overline{\om}}$.

b) For each cocompact lattice $J\subset Z(\B{A})/Z(F)$, 
the quotient $J\bs FBun_{D,n,\overline{\om}}$ a Deligne--Mumford stack, which is a 
quotient of a quasi-projective scheme by a finite group. Furthermore, it is a 
quasi-projective scheme if $J$ is torsion-free, and $|D|$ is sufficiently large 
relative to $\mu$.

c) The induced (from a) actions of  $Z(\B{A})/Z(F)$ on $FBun_{*,n}$ and $FBun_{*,n,\overline{\om}}$
naturally extend to continuous right actions of $\ga/Z(F)$.
\end{Prop}

From now on fix a cocompact lattice $J\subset Z(\B{A})/Z(F)$
(which we may assume to be torsion-free) and an admissible $n$-tuple $\ov{\om}$. 
We are going to define for each $i\in\B{Z}$ the intersection
cohomology with compact support 
$H^i= H^i_J(\ov{\om})=H^i_c(J\bs FBun_{*,n,\overline{\om}}, \IC)$ and its pure quotient
$H_{\pure}^i=H^i_{\pure,J}(\ov{\om})$, both being $\Gm_{F^{(n)}}$-modules over $\ov{\B{Q}_l}$.

\begin{Emp} \label{E:coh}
By \rp{action} b), every $J\bs FBun^{\leq\mu}_{D,n,\overline{\om}}$ is a quotient of a
quasi-projective scheme over $X^n$ by finite group. Therefore we
can consider the intersection cohomology $H_{D}^{i,\mu}:=H_c^i([J\bs
FBun^{\leq\mu}_{D,n,\overline{\om}}]\times_{X^n} F^{(n)}, \IC)$ and
its pure quotient $H_{\pure,D}^{i,\mu}$ (see
\rn{pure}, \ref{N:dm} and \rr{equiv}).

For each $\mu\leq\mu'$ and $D\subset D'$, we have a
natural morphism $H_{D}^{i,\mu}\to H_{D'}^{i,\mu'}$ (using \rp{fbundles} b)), which by Remarks
\ref{R:pure} and \ref{R:equiv} induces an embedding
$H_{\pure,D}^{i,\mu}\hra H_{\pure,D'}^{i,\mu'}$. Thus we can form
direct limits $H^i$ and $H_{\pure}^i$ of the
$H_{D}^{i,\mu}$'s and the $H_{\pure,D}^{i,\mu}$'s, respectively. Both
spaces are equipped with a continuous action of the product $\ga\times \Gm_{F^{(n)}}$.

Finally, let $H_{\pure}^{i,\mu}$ and $H_{\pure,D}^{i}$ be subspaces of $H_{\pure}^i$, obtained as 
direct limits of the $H_{\pure,D}^{i,\mu}$'s taken over $D$'s and $\mu$'s, respectively.
\end{Emp}


\begin{Not} \label{N:orisph}
Let $H^i_{\cusp}=H^i_{\cusp,J}(\overline{\om})$ be the subspace of
$H_{\pure}^i$, consisting of all elements vanishing on all orispheric substacks 
$C\subset J\bs FBun_{*,n,\ov{\om}}$ (recall that orispheric substacks are closed
by \rp{par}, and see \rr{pure} c) for the definition of the
restriction map).
\end{Not}

The main advantage of $H^i_{\cusp}$ over all
previously defined spaces is due to part b) of the following
result, proved in \re{cusp}.

\begin{Prop} \label{P:cusp}
a) For each dominant coweight $\mu\geq d(\overline{\om})\rho$, we have  
\[
H^i_{\cusp,J}(\overline{\om})=\cap_{g\in\ga}g(H_{\pure,J}^{i,\mu}(\overline{\om})).
\]
b) $H^i_{\cusp}$ is an admissible representation of $\ga$.

\end{Prop}

\begin{Conj} \label{Conj}
If $G=GL_r$, then the representation $H^0_{\cusp,J}(\overline{\om})$
of $GL_r(\B{A})\times\Gm_{F^{(n)}}$ is isomorphic to the direct sum
$\bigoplus_{\pi}(\pi\pp\rho_{\pi,\overline{\om}})$, where $\pi$ runs
over the set of all cuspidal representations of $GL_r(\B{A})$ with $\pi(J)=\Id$, and
$\rho_{\pi,\overline{\om}}$ is the same as in the introduction. (Here we identify $\ov{\om}$ with the 
corresponding representation of $(\widehat{G})^n$ as in \rr{corr}.)
\end{Conj}

\begin{Rem} \label{R:cong}
When $G$ is arbitrary, we also expect that cuspidal tempered part of Langlands correspondence  
can be realized in $H^0_{\cusp}$, but one has to take  into account the contribution of 
endoscopic groups as well.
\end{Rem}

To provide an evidence to our conjecture, we will show in Section \ref{S:laf} the following 
result.

\begin{Thm} \label{T:laf}
a) In the Lafforgue case (that is, for $G=GL_r, n=2$ and
$\overline{\om}$ is the tensor product of the standard
representation of $\widehat{G}=GL_r$ and its dual), Conjecture \ref{Conj}
holds up to $r$-negligibles (see \rn{rnegl}). More precisely, there exists an
exhausting filtration $0=V_0\subset V_1\subset
V_2\subset\ldots\subset H^0_{\cusp}$  such that
each $V_{2i+1}/V_{2i}$ is $r$-negligible, and $\bigoplus_{i}
V_{2i}/V_{2i-1}\cong
\bigoplus_{\pi}(\pi\pp\rho_{\pi}\pp\check{\rho}_{\pi})$.

b) Conjecture \ref{Conj} holds in Drinfeld's case (that is, in the Lafforgue case with $r=2$).
\end{Thm}

\begin{Rem} \label{R:gen}
Finally, let us introduce objects, described in the introduction. Let 
$\cal X_n$ and $\cal X_{\overline{\om}}$ be the inverse limits of the $J\bs FBun_{*,n}$'s and 
the  $J\bs FBun_{*,n,\ov{\om}}$'s, respectively, taken over all cocompact lattices in 
$Z(\B{A})/Z(F)$.
Then $H^i_c(\C{X}_{\ov{\om}},\IC)$ and its required subquotient  
$H^i_{\cusp}(\C{X}_{\ov{\om}},\IC)$ are the direct limits of the $H^i_J(\ov{\om})$'s 
and the $H^i_{\cusp,J}(\ov{\om})$'s, respectively. In particular, in the case $G=GL_r$, 
our Conjecture \ref{Conj} for all $J$'s is equivalent to the assertion that 
$H^0_{\cusp}(\C{X}_{\ov{\om}},\IC)$ 
is isomorphic to the direct sum $\bigoplus_{\pi}(\pi\pp\rho_{\pi,\overline{\om}})$, where 
$\pi$ runs over the set of all cuspidal representations with finite order central characters.
\end{Rem}

%

\section{Basic properties of $F$-bundles}
In this  section we will prove Propositions \ref{P:fbundles} and \ref{P:action}. 
For each $\nu\in\pi_1(G)=\pi_0(Bun_G)$, let us denote
by $Bun_{G,D}^{\leq \mu;\nu}$ the preimage in $Bun_{G,D}^{\leq
\mu}$ of the connected component of $Bun_G$,  corresponding to
$\nu$, and  similarly for other spaces such as $Hecke$ and $FBun$. 
We will use the following lemma, whose proof will
be sketched in \re{repr}.

\begin{Lem} \label{L:repr}
a) If $|D|$ sufficiently  large relative to $\mu$, then
$Bun_{G,D}^{\leq \mu;\nu}$ is a smooth quasi-projective scheme for each $\nu\in\pi_1(G)$.

b) Stack $Hecke_{n, \overline{\om}}$ is projective over $Bun\times X^n$.

c) The forgetful map  $\pi:Hecke'_{D,n,\overline{\om}}\to Hecke_{D,n,\overline{\om}}$
is projective. Moreover, $\pi$ is an isomorphism over $X^n\sm\Dt$.
\end{Lem}

\begin{Emp} \label{E:fbundles} 
\begin{proof}[Proof of \rp{fbundles}]
b) The map $(\C{G},\phi)\mapsto (\C{G}_{|D\times S}, \phi_{|D\times S})$ gives a 
morphism from $FBun_{n}\times_{X^n} (X\sm D)^n$ to the stack
$\C{Y}_D$ classifying pairs consisting of a $G$-bundle $\wt{\C{G}}$ on $D\times S$ and an isomorphism 
$^{\tau}\wt{\C{G}}\isom\wt{\C{G}}$. Since $G$ and therefore $G_D$ are geometrically connected,
Lang's theorem implies that $\C{Y}_D$ is a classifying space of the discrete group 
$G_D(\fq)=G(\C{O}_D)$ (use \rl{frob} below). Now the statement follows from the fact that  
$FBun_{D,n}$ is canonically isomorphic to the fiber product of  
$FBun_{n}\times_{X^n} (X\sm D)^n$ and $\Spec \fq$ over $\C{Y}_D\cong G(\C{O}_D)\bs \Spec \fq$.

a) By b), we can replace $D$ by its multiple, so we can assume
that $|D|$ is sufficiently large to satisfy a) of \rl{repr}. Then
$Hecke^{\leq\mu;\nu}_{D,n,\overline{\om}}$, being the restriction of 
$Bun_D^{\leq\mu;\nu}\times_{Bun}Hecke_{n,\overline{\om}}$ to $(X\sm D)^n$, 
is a quasi-projective scheme. Now the statement follows from the
fact that $FBun^{\leq\mu;\nu}_{D,n,\overline{\om}}$ is a closed
substack of $Hecke^{\leq\mu;\nu}_{D,n,\overline{\om}}$. Indeed,
$FBun^{\leq\mu;\nu}_{D,n,\overline{\om}}$ is the preimage in
$Hecke^{\leq\mu;\nu}_{D,n,\overline{\om}}$ of the graph of the
Frobenius morphism in $Bun_{D}^{\leq \mu;\nu}\times
Bun_{D}^{\leq \mu;\nu}$.

c) We have to check that $Bun_{D}$ satisfies the conclusion of \rl{frob} b) below.
Instead of checking that $Bun_D$ satisfies the assumption of  \rl{frob} b), we can argue as follows. 
As the question is local for the Zariski topology,
we may replace $Bun_D$ by its open substack $Bun_D^{\leq\mu;\nu}$. By \rl{repr} a), 
there exists a finite subscheme $D'\subset X$ containing $D$ such that  
$Bun_{D'}^{\leq\mu;\nu}$ is a scheme. Then 
$Bun_{D'}^{\leq\mu;\nu}$ clearly satisfies \rl{frob}, so the statement for  
$Bun_D^{\leq\mu;\nu}$ follows from b) together with the fact that groupoid 
$Bun^{\leq\mu;\nu}_{D}(\fq)$ is isomorphic to the quotient of 
$Bun^{\leq\mu;\nu}_{D'}(\fq)$ by $G_{D,D'}(\fq)$ (use again Lang's theorem).


d) The ``only if'' statement was explained in \rr{adm}.
Assume now that $\overline{\om}$ is admissible. Then $FBun_{D,n,\overline{\om}}$ 
contains a substack consisting of $F$-bundles for which
$x_1=\ldots=x_n$ and $\phi$ is an isomorphism. As this substack
obviously contains (and actually is isomorphic by c) to)
a non-empty stack $Bun_D(\fq)\times (X\sm D)$, the statement follows.

e) follows immediately from statement c) of  \rl{repr}.
\end{proof}
\end{Emp}

\begin{Lem} \label{L:frob}
a) Let $\C{X}$ be an algebraic stack locally of finite type over $\fq$, and let $\C{Y}$
be a stack over $\fq$ such that 
$\C{Y}(S)=\{(A,\Phi)\,| A\in\C{X}(S), \Phi\in Isom_{\C{X}(S)}({}^{\tau}A,A)\}$. 
Then  $\C{Y}$ is a Deligne--Mumford stack, \'etale over $\fq$, containing the discrete stack 
 $\C{X}(\fq)$ as an open and closed substack.

b) If, in the notation of a), all geometric fibers of the diagonal morphism \break
$\Dt_{\C{X}}:\C{X}\to\C{X}\times\C{X}$ are connected, then $\C{Y}$ is canonically isomorphic
to $\C{X}(\fq)$.
\end{Lem}
\begin{pf}
a) As $\C{Y}$ is a fiber product  over $\C{X}\times\C{X}$ of the diagonal $\Dt_{\C{X}}$ and the 
graph of Frobenius morphism, $\C{Y}$ is an algebraic stack locally of finite type over $\fq$. 
One checks that the canonical morphism $\C{Y}(\ov{\fq}[t]/(t^2))\to\C{Y}(\ov{\fq})$ 
is an equivalence of categories, therefore the diagonal morphism $\Dt_{\C{Y}}$ is unramified. 
Hence $\C{Y}$ is a Deligne--Mumford stack  (see \cite[Thm. 8.1]{LMB}), 
and $\C{Y}$ is \'etale over $\fq$. 

It remains to check that the natural functor 
$i:\C{X}(\fq)\to\C{Y}(\ov{\fq})$ is fully faithful. But this follows from (actually is equivalent to) 
the first axiom of a stack (sheaf axiom for $\Isom(x,y)$)
applied to \'etale covers $\Spec \B{F}_{q^m}\to \Spec \fq$ for all $m\in\B{N}$. 

b) We have to show that the functor $i:\C{X}(\fq)\to\C{Y}(\ov{\fq})$ from a) is essentially surjective. 
Let $(A,\Phi)$ be any object of
$ \C{Y}(\ov{\fq})$, and we want to find an object $B$ of $\C{X}(\fq)$ such that $i(B)$ is 
isomorphic to $(A,\Phi)$. Choose $m$ such that $A$ is (a pull-back of) an object of 
$\C{X}(\B{F}_{q^m})$, and $\Phi$ belongs to $\Isom_{\C{X}(\fqm)}({}^{\tau}A,A)$.
Then 
\[
\Phi^{(m)}:=\Phi\circ{}^{\tau}\Phi\circ\ldots\circ{}^{\tau^{m-1}}\Phi:A={}^{\tau^m}A\isom A
\]
defines an $\fqm$-point of an algebraic group $H:=\Isom_{\C{X}}(A,A)$. 
If $\Phi^{(m)}$ is the identity, then the existence of $B$ is equivalent to the second
axiom of a stack applied to the \'etale cover  $\Spec \B{F}_{q^m}\to\Spec \fq$.
The general case easily reduces to this one. Indeed, our assumption about 
$\Dt_{\C{X}}$ implies that $H$ is a connected group over $\B{F}_{q^m}$. 
Therefore Lang's theorem implies the existence of $h\in H(\ov{\fq})$
such that $^{\tau^m}h=h\circ\Phi^{(m)}$. Then $h$ induces an isomorphism between 
$(A,\Phi)$ and $(A,\Phi':=h\circ\Phi\circ {}^{\tau}\!h^{-1})$. By construction,  
$\Phi'^{(m)}$ is the identity, completing the reduction.
\end{pf}

\begin{Rem} \label{R:frob}
a) \rl{frob} gives an explanation of the Drinfeld's lemma 
(see for example \cite[Ch.1, 3, Lem. 3]{La1}) used by Drinfeld and Lafforgue.

b) In the proof of \rp{fbundles} we used \rl{frob} only in two particular cases: 
when $\C{X}$ is a scheme and when $\C{X}$ is a classifying space of a connected algebraic group.
In both cases the proof can be simplified.
\end{Rem}

\begin{Emp} \label{E:action1}
Before beginning the proof of \rp{action}, recall that $G(\B{A})$ acts naturally (on the right) on 
the inverse limit $Bun_{*}$ of the $Bun_{D}$'s. Any element of $Bun_{*}(S)$ consists of a
$G$-bundle $\C{G}$ over $X\times S$ equipped with trivializations
$\phi_v:\C{G}_{|\C{O}_v\times S}\isom G\times\C{O}_v\times S$ for
all closed points $v$ of $X$. For each $g=(g_v)_v\in G(\B{A})$,
choose a finite set of closed points $T$ of $X$ such that $g_v\in
G(\C{O}_v)$ for all $v\notin T$. We claim that there exists a
unique $G$-bundle $\wt{\C{G}}$ over $X\times S$ such that
$\wt{\C{G}}_{|(X\sm T)\times S}=\C{G}_{|(X\sm T)\times S}$ and
$\wt{\C{G}}_{|\C{O}_v\times
S}=\phi_v^{-1}(g_v(G\times\C{O}_v\times S))$ for each $v\in T$
(the last equality we consider inside 
$\wt{\C{G}}_{|F_v\times S}=\C{G}_{|F_v\times S}$). Indeed, by
\cite{BL}, the corresponding statement holds for vector bundles,
so by Tannakian formalism it holds in general. 

Since $\wt{\C{G}}$ is clearly independent of $T$, the rule
$(\C{G},\{\phi_v\}_v)g:=(\wt{\C{G}}, \{g_v^{-1}\circ\phi_v\}_v)$
defines the required group action. Moreover, it follows from the construction that 
$Z(F)$ acts trivially, and that the induced action
of $Z(\B{A})/Z(F)$ on $Bun_{*}$ gives an action on each $Bun_{D}$. Furthermore,
since the open substack $Bun_G^{\leq\mu}$ was defined as the preimage of 
$Bun_{G^{\ad}}^{\leq\mu}\subset Bun_{G^{\ad}}$, the group $Z(\B{A})/Z(F)$ preserves each 
$Bun_G^{\leq\mu}$.
\end{Emp}

\begin{Rem} \label{R:action}
The above argument actually shows that $Bun_*$ is equipped with an action of a huge 
ind-pro-algebraic group, whose group of $\fq$-points is $G(\B{A})$.
\end{Rem}

\begin{Emp} \label{E:action}
\begin{proof}[Proof of \rp{action}]
a) and c) follow from the fact that the action of $G(\B{A})/Z(F)$ on $Bun_*$ (resp.
 $Z(\B{A})/Z(F)$ on $Bun^{\leq\mu}_D$), defined in \re{action1}, naturally lifts to the 
actions on $FBun_{*,n}$ (resp.  $FBun^{\leq\mu}_{D,n}$) and leaves  $FBun_{*,n,\ov{\om}}$ 
(resp. $FBun^{\leq\mu}_{D,n,\ov{\om}}$) invariant.

b) By \rp{fbundles} b), the statement would follow if we show that the induced action of $J$ 
on $\pi_0(FBun_{D,n,\ov{\om}})$ has finitely many orbits and has finite stabilizers.
As we proved in \re{fbundles} that the projection 
$\pi_0(FBun_{D,n,\overline{\om}})\to\pi_0(Bun)$
has finite fibers, it will suffice to prove the corresponding statement for $\pi_0(Bun_G)$
instead of $\pi_0(FBun_{D,n,\ov{\om}})$. 

Notice that $Z(\B{A})$ acts on $\pi_0(Bun_G)=\pi_1(G)$ via the continuous 
homomorphism $\Pi:Z(\B{A})\hra T(\B{A})\overset{\pi}{\longrightarrow}
X_*(T)\to\pi_1(G)$, whether $\pi$ is given by the rule
$\langle \pi(t),\nu\rangle =\log_q|\nu(t)|$ for each $\nu\in X^*(T)$ and $t\in
T(\B{A})$. Indeed, as the natural surjection $\pi_0(Bun_T)=X_*(T)\to\pi_1(G)=\pi_0(Bun_G)$ is
induced by the inclusion $T\hra G$ (see the proof of \rl{conn} in \re{conn}), 
it will suffice to show the corresponding statement for $G=T$, hence for $G=\B{G}_m$, in which
case it is clear.

Note that $\Pi$ factors through $Z(\B{A})/Z(F)$, and that the induced homomorphism 
$\Pi':Z(\B{A})/Z(F)\to\pi_1(G)$ has a compact kernel and a finite cokernel.
Since $J\subset Z(\B{A})/Z(F)$ is a cocompact lattice, both kernel and cokernel
of the restriction of $\Pi'$ to $J$ are therefore compact and discrete, hence finite.
This implies the assertion.
\end{proof}
\end{Emp}


\section{Local model of $FBun_{n,\ov{\om}}$}

The goal of this section is to prove \rt{locisom} and \rco{strata}. 
Our strategy will be to decompose locally
$Hecke_{n,\overline{\om}}$ as a product
$Gr_{n,\overline{\om}}\times Bun$  and then to use the fact that
the Frobenius morphism has a zero differential.

\begin{Lem} \label{L:locisom}

Let  $\C{G}_0$ be a  $G$-bundle on $X\times S$, locally trivial in
the Zariski topology, and let $\pi:S\to Bun_{G}$ be the morphism,
corresponding to $\C{G}_0$. 

Then the fiber product
$Hecke_{n,\overline{\om}}\times_{Bun} S$ 
(taken with respect to the projection $p'$ (see \rn{proj}))
 and the product
$Gr_{n,\overline{\om}}\times S$ are locally Zariski isomorphic
fibrations over $X^n\times S$. Moreover, the isomorphism preserves
the stratifications induced by those of $Hecke_{n,\overline{\om}}$
and $Gr_{n,\overline{\om}}$. 

Furthermore, this local isomorphism lifts to a stratification preserving local isomorphism between  
$Hecke'_{n,\overline{\om}}\times_{Bun} S$ and
$Gr'_{n,\overline{\om}}\times S$.
\end{Lem}

\begin{pf}
Let  $(x'_1,\ldots,x'_n;s')$ be a closed point of $X^n\times S$. We
want to find its open neighborhood, whose  inverse images in
$Hecke_{n,\overline{\om}}\times_{Bun} S$ and
$Gr_{n,\overline{\om}}\times S$ are isomorphic. 
We are going to prove the statement by induction on $n$.

Assume first that $x'_1=x'_2=\ldots=x'_n$ (this condition holds automatically
for $n=1$).
By our assumption, there exists an open neighborhood
$V\subset X\times S$ of $(x'_1,s')$ and a trivialization
$\psi$ of the restriction of $\C{G}_0$ over $V$. 
Consider the open subscheme $U$ of $X^n\times S$, consisting of points
$(x_1,\ldots,x_n;s)$ such that $(x_i,s)\in V$ for each $i=1,\ldots,n$.
We claim that $U$ is a required neighborhood. Let
$U'\subset Hecke_{n,\overline{\om}}\times_{Bun} S$ and $U''\subset
Gr_{n,\overline{\om}}\times S$ be the  inverse images of $U$. We
are going to find an isomorphism $U'\isom U''$ over $U$.

Let $(\C{G},\C{G}'; y_1,\ldots,y_n;\phi)$ be the pullback to $U'$ of
the universal object over $Hecke_{n,\overline{\om}}$. Let
$V'\subset X\times U'$ be the preimage of $V$, and 
let $\psi'$ be the trivialization of $\C{G}'$ over
$V'$ induced by $\psi$. The composition of $\phi$ and
$\psi'$ defines a trivialization $\phi'$ of $\C{G}$ over
$V'\sm (\Gm_{y_1}\cup\ldots\cup \Gm_{y_n})$.
As $V'$ contains each $\Gm_{y_i}\subset X\times U'$,
there exists a unique $G$-bundle $\wt{\C{G}}$ over
$X\times U'$, trivial over $(X\times U')\sm (\Gm_{y_1}\cup\ldots\cup
\Gm_{y_n})$ such that $\wt{\C{G}}_{|V'}=\C{G}_{|V'}$ and the gluing is done
by means of $\phi'$.
By the definition of the affine grassmannian,
$\wt{\C{G}}$ defines the required morphism $U'\to U''$. The
construction of the inverse map and rest of the statements is now straightforward and 
is therefore omitted (compare the proof of \rl{red}).  

Assume now that $x'_i\neq x'_j$ for some $i$ and $j$; then after renumbering 
of indexes there exists a positive integer $k<n$ such that $x'_i\neq x'_j$ for  
each $i\leq k<j$. Using the second statement of \rl{red} a), the assertion now follows by induction hypothesis.
\end{pf}

\begin{Emp} \label{E:locmod}
\begin{proof}[Proof of \rt{locisom}]
Choose $\mu$ such  that $y\in FBun^{\leq\mu}_{D,n,\overline{\om}}$.
By \rp{fbundles} b) and \rl{repr} a), we can enlarge $D$ (and
replace  $y$ by one of its preimages) so that $Bun^{\leq\mu}_{D}$ is a scheme. By the
theorem of Drinfeld--Simpson (\cite{DS}), there exists a surjective
\'etale morphism $\pi: S\to Bun^{\leq\mu}_D$ such that the pull-back
to $X\times S$ of the universal $G$-bundle on $X\times
Bun_D$ is locally trivial in the Zariski topology.

Choose  a preimage  $y'\in FBun_{D,n,\overline{\om}}\times_{Bun_D}
S$ of $y$, and let $y''\in X^n\times S$ be the
image of $y'$. By \rl{locisom}, $y''$ has an open neighborhood $U\subset X^n\times S$,
whose inverse images  $U'\subset Hecke_{n,\overline{\om}}\times_{Bun} S$  
and $U''\subset Gr_{n,\overline{\om}}\times S$ are isomorphic over $U$.
We claim that $U_y:=U'\times_{Hecke_{D,n}} FBun_{D,n}$ is the required \'etale neighborhood
of $y$. Since the natural projection $U_y\to  FBun_{D,n}$ is \'etale (since
$\pi$ is so), it will suffice to show that the 
composition map $U_y\hra U'\isom U''\hra Gr_{n,\overline{\om}}\times
S\to  Gr_{n,\overline{\om}}$ is \'etale.

Denote by $f$ the restriction to $U'$ of the projection $p:Hecke_{D,n}\to Bun_D$
from \rn{proj}. Then if we identify $U'$ and $U''$ by means of isomorphism 
$U'\isom U''$, chosen above, the statement follows from \rl{etale} below, applied to 
$Y:=S$, $Z:=Bun^{\leq\mu}_{D}$, $T:= Gr_{n,\overline{\om}}$ and $W:=f^{-1}(Bun^{\leq\mu}_{D})$ 
(hence $V=U_y$).

Furthermore, by the last assertion of \rl{locisom},  $U_y\to Gr_{n,\overline{\om}}$ lifts to 
a morphism $U_y\times_ {FBun_{D,n,\overline{\om}}} FBun'_{D,n,\overline{\om}}\to Gr'_{n,\overline{\om}}$,
which is again \'etale by \rl{etale} below.
\end{proof}
\end{Emp}

\begin{Lem} \label{L:etale}
Let $Y,Z$ and $T$ be schemes locally of finite type over $\fq$, let $W\subset Y\times T$
be an open subscheme, let $\pi:Y\to Z$ be an \'etale morphism, let $f:W\to Z$ be any morphism, and
let $V$ be given by equation 
$$V=\{(y,t)\in W\subset Y\times T| Frob_q(f(y,t))=\pi(y)\}.$$
Assume that $Z$ is smooth over $\fq$. Then the canonical map $\Pi:V\hra W\to T$ is \'etale.
\end{Lem}

\begin{proof}
Assume first that $T$ is smooth over $\fq$ (compare \cite[Prop. 3.3]{Dr1}).
In this case $Y$ and $W$ are smooth as well. Since the Frobenius morphism has a zero differential, 
$V$ inside $W$ is locally given by $\dim Z$ equations with linearly independent differentials. 
Therefore $V$ is smooth. Moreover, the projection $\Pi:V\to T$ induces an isomorphism on tangent 
spaces, hence $\Pi$ is \'etale, as claimed. 

In the general case, we may assume that $Z=\B{A}^m$. Indeed, 
as the question is local on $V$, we may shrink all the schemes in question so that  all of them are affine 
and there exists an \'etale map $\phi:Z\to\B{A}^m$. Then $V$ is an open and closed 
subscheme of the scheme
$$V'=\{(y,t)\in W\subset Y\times T| Frob_q(\phi\circ f(y,t)))=\phi\circ\pi(y)\}.$$ 
In particular, we may replace $Z$, $\phi$ and $f$ by $\B{A}^m$, $\phi\circ\pi$ and $\phi\circ f$, 
respectively.

Next choose a closed embedding of $T$ into an affine space $\wt{T}=\B{A}^k$. Then there exists 
an open affine subscheme 
$\wt{W}\subset Y\times \wt{T}$ such that $\wt{W}\cap (Y\times T)=W$ and an extension 
$\wt{f}:\wt{W}\to Z=\B{A}^m$ of $f$ to $\wt{W}$. Consider 
$$
\wt{V}:=\{(y,t)\in \wt{W}\subset Y\times \wt{T}| Frob_q(\wt{f}(y,t))=\pi(y)\}.
$$
As $\wt{T}$ is smooth, we have seen before that the natural projection $\wt{\Pi}:\wt{V}\to\wt{T}$ is \'etale. 
But $\Pi:V\to T$ is just the restriction of $\wt{\Pi}$ to $V=\wt{\Pi}^{-1}(T)$. 
Therefore it is also \'etale, as claimed.
\end{proof}

\begin{Emp} \label{E:strata}
\begin{proof}[Proof of \rco{strata}]

b) and the first statement of a) are reduced by \rt{locisom} to the corresponding questions about 
affine grassmannians, which will be shown in  \rp{affgr}. As for the second statement of a), the ``only if'' part
was explained in \rr{fbundles} b), while the ``if'' part  
for  $FBun^0_{D,n,\overline{\om}}$ follows  from \rp{fbundles} d). Moreover, since 
$\pi:FBun'_{D,n,\overline{\om}}\to FBun_{D,n,\overline{\om}}$ induces an isomorphism 
$FBun'^0_{D,n,\overline{\om}}\times_{X^n} (X^n\sm\Dt)\isom FBun^0_{D,n,\overline{\om}}$ (by \rp{fbundles} e)), 
the ``if'' part for $FBun'^0_{D,n,\overline{\om}}$ follows as well.

c) Let $m$ be the dimension of $Bun_D$, and let $\C{F}(m/2)[m]$ 
be the restriction of the IC-sheaf of
$Hecke_{D,n,\overline{\om}}$ to $FBun_{D,n,\overline{\om}}$. We
want to show that $\cal F$ is the IC-sheaf. Since the statement clearly holds for 
the restriction of $\C{F}$ to the smooth open dense stratum $FBun^0_{D,n,\overline{\om}}$,
it remains to show that $\C{F}$ is an irreducible perverse sheaf.

As the map from the disjoint union of the $U_y$'s (from  \rt{locisom}) to
$FBun_{D,n,\overline{\om}}$ is \'etale and surjective, it will suffice to show the corresponding
statement for the restriction of 
$\cal F$ to each $U_y$. Consider the commutative diagram, 

$$
\CD
U_y              @>{i}   >> U'@>{\pi_1}>>Gr_{n,\overline{\om}}\\
@VVV                  @V{\pi_2}VV\\
FBun_{D,n,\overline{\om}}@>>> Hecke_{D,n,\overline{\om}}
\endCD
$$
constructed in the course of the proof of \rt{locisom}.
 As  $\pi_2$ is \'etale, $\pi_1$ is smooth of relative dimension $m$, and $\pi_1\circ i$ is \'etale, 
we get that 
$$
\C{F}_{|U_y}=i^*\pi_2^*(\IC_{Hecke_{D,n,\overline{\om}}}(-\frac{m}{2})[-m])=
i^*(\IC_{U'}(-\frac{m}{2})[-m])=i^*\pi_1^*(\IC_{Gr_{n,\ov{\om}}})=\IC_{U_y},
$$ 
as claimed.

The last assertion follows immediately from the corresponding 
statement for $Hecke_{n,\overline{\om}}$ shown in \rp{tate} and the observation that
the difference \break
$\dim FBun_{D,n,\overline{\om}}-n =
2\langle \sum^n_{i=1} \om_i,\rho\rangle$ is even.

Finally, the assertion for $FBun'_{D,n,\overline{\om}}$ follows from the (proof of) the corresponding 
assertion of \rt{locisom} by precisely the same argument.

d) By \rp{fbundles} e), $\pi:FBun'_{D,n,\overline{\om}}\to FBun_{D,n,\overline{\om}}$ is projective and 
induces an isomorphism 
$FBun'^0_{D,n,\overline{\om}}\times_{X^n} (X^n\sm\Dt)\isom FBun^0_{D,n,\overline{\om}}$. Hence the 
surjectivity statement follows from a). Finally, by \rt{locisom} and \rl{locisom}, 
the smallness of $\pi$ is equivalent to the smallness of the forgetful morphism
$Hecke'_{n,\overline{\om}}\to Hecke_{n,\overline{\om}}$, which will be shown in \rl{small}.
\end{proof}
\end{Emp}

\section{Reducible $F$-bundles}

In order to prove \rt{red}, we first need some preparations.
Fix a $G$-bundle $\C{G}$ over a smooth connected
projective curve $X$ over an algebraically closed field $k$. 
Equip the set of all $B$-structures of $\C{G}$ with a following partial order: 
we say that $\C{B}'\leq\C{B}''$ if $\deg(\C{B}'_{\la})\leq \deg(\C{B}''_{\la})$
for each dominant (or equivalently quasi-fundamental) weight $\la$ of $G$.

\begin{Lem} \label{L:uniq}
a) For every $B$-structure $\C{B}$  of $\C{G}$ there exists a $B$-structure $\C{B}'$ such that
$\C{B}\leq\C{B}'$ and $\C{B}'$ is maximal  with respect to the
above order.

b) Let $\C{B}$ be a $B$-structure of $\C{G}$, which is maximal with respect to the
above order, and let 
$\C{P}$ be any maximal parabolic structure of $\C{G}$, corresponding to a certain simple root 
$\al$ of $G$ with the corresponding quasi-fundamental weight $\la$. 
Then either $\C{P}$ contains $\C{B}$, 
or $\deg(\C{P}_{\la})\leq \deg(\C{B}_{\la})-\deg(\C{B}_{\al})\langle \la,\al\rangle +
4g\langle \rho,\la\rangle $.
\end{Lem}

\begin{Rem} \label{R:constant}
The inequality is very far from being optimal.
\end{Rem}

%
%
%
%

Let us first show how \rl{uniq} implies \rt{red}.

\begin{Emp} \label{E:red}
\begin{proof}[Proof \rt{red}]
Let $(\C{G};x_1,\ldots,x_n;\phi)$ be any geometric point of the complement 
$FBun_{n,\ov{\om}}\sm FBun_{n,\ov{\om}}^{\leq d(\ov{\om})\rho}$.  
Then there exists a $B$-structure $\C{B}$ of $\C{G}$ and a dominant weight $\la'$ of $G$ 
such that $\deg(\C{B}_{\la'})>\langle d(\overline{\om})\rho,\la'\rangle $. By \rl{uniq} a), 
we may assume that $\C{B}$ is maximal with respect to the above order.
Hence there exists a simple root $\al$ of $G$ such that
$\deg(\C{B}_{\al})>\langle d(\overline{\om})\rho,\al\rangle =d(\overline{\om})$. 
Let $\C{P}\supset\C{B}$ be the maximal parabolic structure of $\C{G}$, corresponding to $\al$. 
We want to show that $\phi$ induces a rational
isomorphism between $\C{P}$ and ${}^{\tau}\C{P}$. Clearly, $\phi^{-1}$ induces a rational 
isomorphism between  ${}^{\tau}\C{P}$ and a certain parabolic structure $\C{P}'$ of $\C{G}$,
so it remains to show that $\C{P}'$ contains $\C{B}$.

Let $\la$ be the quasi-fundamental weight of $G$, 
corresponding to $\al$. By \rl{par}, $\C{P}$ and $\C{P}'$ define line subbundles
$\C{L}$ and $\C{L}'$ of $\C{G}_{\la}$, respectively, satisfying Pl\"ucker relations. As
$\C{L}\cong\C{B}_{\la}$ and $\C{L}'\cong\C{P}'_{\la}$, the statement will follow from 
\rl{uniq} applied to $\C{P}'$ if we check that 
$\deg(\C{L}')\geq \deg(\C{L})-d(\overline{\om})\langle \la,\al\rangle +4g\langle \rho,\la\rangle $.

By \rr{inv} c), $\phi^{-1}({}^{\tau}\C{G}_{\la})$ is contained in
$\C{G}_{\la}(\sum_{k=1}^n\langle \om_k,-w_0({\la})\rangle x_i)$. Therefore
$\phi^{-1}({}^{\tau}\C{L})\subset\C{L}'
(\sum_{k=1}^n\langle \om_k,-w_0({{\la}})\rangle x_i)$.
It follows that  
\[
\deg(\C{L}')\geq \deg(\C{L})-\langle \sum_{k=1}^n\om_k,-w_0({\la})\rangle;
\]
thus it remains to check that $\langle \sum_{k=1}^n \om_k,-w_0(\la)\rangle 
\leq d(\overline{\om})\langle \la,\al\rangle -4g\langle \rho,\la\rangle $. 
Since  $-w_0(\la)$ is the quasi-fundamental weight of $G$ corresponding to the
simple root $-w_0(\al)$, the statement follows from the definition of $d(\ov{\om})$ 
and the equality $\langle \rho,\la\rangle =\langle -w_0(\rho),\la\rangle =
\langle \rho,-w_0(\la)\rangle $.
\end{proof}
\end{Emp}

\begin{Emp} \label{E:uniq} 
\begin{proof}[Proof of \rl{uniq}] 
a) It will suffice to show that for each quasi-fundamental weight $\la$ of $G$, the set
 $\{\deg(\C{B}_{\la})\}_{\C{B}}$ is bounded from above. Since every $\C{B}_{\la}$ is 
 canonically a line subbundle of $\C{G}_{\la}$, the statement follows.

b) Let $\C{B}$ be any maximal $B$-structure of $\C{G}$. First we claim that 
$\deg(\C{B}_{\beta})\geq -2g$ for every simple root $\beta$ of $G$. Assume first that 
$G=GL_2$. In this case, our
claim asserts that every rank two vector bundle $\C{E}$ contains a line subbundle of 
degree at least $\frac{1}{2}\deg(\C{E})-g$, so it follows immediately from the Riemann--Roch 
theorem. The general case reduces to that of $GL_2$. Indeed, fix any $\beta$. Let
$P_{\beta}\supset B$ be the  parabolic subgroup $G$ such that $\beta$ is the only simple root 
of its Levi subgroup, and let $R(P_{\beta})$ be the radical of $P_{\beta}$.  
Consider Borel structure $\C{B}':=R(P_{\beta})\bs\C{B}$ of the $P_{\beta}/R(P_{\beta})$-bundle  
$R(P_{\beta})\bs [P_{\beta}\times_B\C{B}]$. Then $\deg(\C{B}'_{\beta})=\deg(\C{B}_{\beta})$, and 
$\C{B}'$ is maximal. [If not, then there exists a Borel structure $\C{B}'_0$ larger 
than $\C{B}'$. Hence the preimage of  $\C{B}'_0$ in $P_{\beta}\times_B\C{B}$ would give us a 
$B$-structure of $\C{G}$, larger than $\C{B}$, contradicting the maximality of $\C{B}$].
Since $P_{\beta}/R(P_{\beta})\cong PGL_2$, we are thus reduced to the case of $PGL_2$, hence to 
that of $GL_2$, as claimed.

Now we are ready to prove the assertion. Observe first that we can replace $G$
by $G^{\ssc}$ and $\la$ by the corresponding fundamental weight. 
By \rl{par}, $\C{B}$ and $\C{P}$ define line 
subbundles $\C{L}$ and $\C{L}'$ of $\C{G}_{\la}$, respectively, 
and we have to show that either $\C{L}'=\C{L}$
or $\deg(\C{L}')\leq \deg(\C{L})-\deg(\C{B}_{\al})+ 4g\langle \rho,\la\rangle $. 
For this we will show that the latter inequality holds for every
line subbundle $\C{L}'\neq\C{L}$ of  $\C{G}_{\la}$. Fix any $B$-invariant complete flag 
$0=V_0\subset V_1\subset\ldots\subset V_M=V_{\la}$ of $V_{\la}$.
Then $\C{B}$ defines a complete flag 
$0=\C{E}_0\subset\C{E}_1\ldots\subset\C{E}_M=\C{G}_{\la}$ of $\C{G}_{\la}$ with 
$\C{E}_i=B\bs[\C{B}\times V_i]$. By construction, each quotient $\C{E}_i/\C{E}_{i-1}$ is isomorphic to 
$\C{B}_{\la_i}$, where $\la_i$ is the weight of $V_i/V_{i-1}$. In particular, $\C{E}_1\cong\C{L}$.

For each line subbundle $\C{L}'\neq\C{L}$, let
$i\geq 2$ be the smallest integer such that $\C{L}'$ is contained in $\C{E}_i$. Then $\C{L}'$ 
embeds into some $\C{E}_i/\C{E}_{i-1}\cong\C{B}_{\la_i}$, so it will suffice to show that 
$\deg(\C{B}_{\la'})\leq \deg(\C{B}_{\la})-\deg(\C{B}_{\al})+
4g\langle \rho,\la\rangle $ for every weight $\la'\neq\la$ of $V_{\la}$ or, equivalently, that 
$\deg(\C{B}_{\la-\la'-\al})\geq -4g\langle \rho,\la\rangle$. But  $\la-\la'-\al$ 
is of the form $\sum_{\beta}n_{\beta}\beta$, where $\beta$ runs over the set of all 
simple roots of $G$ and each $n_{\beta}$ is non-negative.
Since $\deg(\C{B}_{\beta})\geq -2g$ for each $\beta$, it remains to show that 
$\sum_{\beta}n_{\beta}\leq 2\langle \rho,\la\rangle$. As
$$
\sum_{\beta}n_{\beta}=\langle \la-\la'-\al,\rho\rangle< \langle \la-\la',\rho\rangle
\leq \langle \la-w_0(\la),\rho\rangle =\langle \la,\rho-w_0(\rho)\rangle =2\langle \rho,\la\rangle, 
$$ 
we get the assertion.
\end{proof}

\end{Emp}

To formulate a more precise version of \rp{par}, we will need the following result.

\begin{Cl} \label{C:dec}
Fiber product $FBun_{P,D,n}\times _{X^n}(X^n\sm\Dt)$ naturally decomposes as a disjoint union of 
open and closed substacks  $FBun_{P,D,n}^{d;\bar{k};[g]}$, 
indexed by triples $(d;\bar{k};[g])$, where $d$
is an integer, $\bar{k}$ is an $n$-tuple of integers with zero sum,
and $[g]$ is an element of $G(\C{O}_D)/P(\C{O}_D)$. 
\end{Cl}

\begin{proof}
Let $\C{G}$ be the universal $G$-bundle on  $X\times FBun_{P,n}$,
and let $\la$ be the quasi-fundamental weight of $G$ corresponding
to $P$. By \rl{par}, the universal $P$-bundle on $FBun_{P,n}$
defines a line subbundle $\C{L}\subset\C{G}_{\la}$, satisfying
Pl\"ucker relations and such that the rational isomorphism
$\phi_{\la}$ between $\C{G}_{\la}$ and ${}^{\tau}\C{G}_{\la}$
induces that between  $\C{L}$ and ${}^{\tau}\C{L}$. 

Now we claim that for each pair $(d,\bar{k})$ as in the assertion,
there exist an open and closed substack
$FBun_{P,n}^{d;\bar{k}}$ of  $FBun_{P,n}\times_{X^n}(X^n\sm\Dt)$
consisting of geometric points $s$ such that $\deg(\C{L}_s)=d$ and
${}^{\tau}\C{L}_s=\phi(\C{L}_s)(\sum_{i=1}^n k_i x_i)$.
The first condition is clearly open and closed. For the second, 
consider the unique line bundle $\C{L}_i$ on $X\times [FBun_{P,n}\times_{X^n}(X^n\sm\Dt)]$, 
whose restriction to the complement of $\Gm_{x_i}$ is
${}^{\tau}\C{L}$ and whose restriction to the complement of  
$\cup_{j\neq i}\Gm_{x_j}$ is $\phi(\C{L})$. Since  
the second condition is equivalent to $\deg(\C{L}_i)_s-\deg(\C{L}_s)=k_i$ for each $i$,
it is open and closed as well. 

Finally, observe that the map $(\C{G},\C{P},\psi,...)\mapsto \psi(\C{P}_{|D\times S})$
defines a morphism from  $FBun_{P,D,n}$ to the stack classifying Frobenius-equivalent
$P_D$-structures of the trivial $G_D$-structure. As the latter stack is isomorphic to the 
discrete stack $G(\C{O}_D)/P(\C{O}_D)$ (compare \rl{frob}), we
get the required decomposition by the $[g]$'s.
\end{proof}

\begin{Rem}
Passing to the limit over $D$'s, we get a decomposition of
$FBun_{P,*,n}$ indexed by the triples as above but with $[g]$'s 
belonging to $G(\B{O})/P(\B{O})$, where $\B{O}$ is the ring of integral
adeles of $F$.
\end{Rem}


\begin{Prop} \label{P:par1}
For each triple $(d;\bar{k};[g])$ as in \rcl{dec}, 
the restriction of the projection $\Pi:FBun_{P,D,n}^{d;\bar{k};[g]}\to FBun_{D,n}$ to the generic fiber
over $X^n$ is finite and unramified. Furthermore, each projection
$FBun_{P,*,n}^{d;\bar{k};[g]}\to FBun_{*,n}$ is a closed
embedding. In particular, every orispheric substack of
$FBun_{*,n}$ is closed.
\end{Prop}

\begin{proof}

As $\G(\B{O})$ acts transitively on the set of $[g]$'s, we may and will assume that
$[g]=[1]$. Since Pl\"ucker relations are closed, \rl{par} implies that it is
enough to show the statement in the case $G=GL_m$ and $P$ is the
maximal parabolic corresponding to the standard representation. 
Also the statement is clearly local on the base, and the first assertion 
is independent of $D$ (by \rp{fbundles} b)). 
Thus it will suffice to check that for each quasi-compact open substack $V$ of $Bun_{GL_m}$, 
each sufficiently large $D$ and each sufficiently small open subscheme $U\subset X^n$ 
(both depending on $V$), the restriction of $\Pi$ to the preimage of 
$V\times U\subset Bun_{GL_m}\times X^n$ is a closed embedding.

Given $V$, let $l_1$ and $l_2$ be two integers such that for every geometric point
of $V$, the corresponding vector bundle does not have line (resp. rank two) subbundles of degree greater 
than $l_1$ (resp. $l_2$). Let $D\subset X$
be a finite subscheme such that $|D|>\max\{l_1-d,l_2-2d\}$ and $D$ contains all points of $X$ of
degree $\leq(l_1-d)$ over $\fq$. Finally let $U\subset X^n$ be an open subscheme such that each 
$(x_1,\ldots,x_n)\in U$ satisfies $x_i\neq {}^{\tau^r}x_j$ for each $i,j$ and each
$r=1,\ldots,l_1-d$. Denote the preimages of $V\times U\subset Bun_{GL_m}\times X^n $ 
in  $FBun_{D,n}$ and $FBun_{P,D,n}^{d;\bar{k};[g]}$ by $\C{A}$ and $\C{B}$, respectively,
and we going to check that the projection $\C{B}\to\C{A}$ is a closed embedding.

Consider the natural morphism $\nu:\C{B}\to FBun^{d}_{GL_1,D,n,\bar{k}}\times_{(X\sm D)^n}\C{A}$,
where the first projection $\C{B}\to FBun^{d}_{GL_1,D,n,\bar{k}}$ was defined during the proof of
\rcl{dec}. Note that $FBun^{d}_{GL_1,D,n,\bar{k}}\times_{(X\sm D)^n}\C{A}$ classifies pairs consisting of 
a line bundle $\C{L}$ and a rank $m$ vector bundle $\C{E}$ on $X\times S$, equipped with $D$-level
structures and $F$-structures (that is, rational isomorphisms from $\C{L}$ and $\C{E}$ to ${}^{\tau}\C{L}$
 and ${}^{\tau}\C{E}$, respectively). 
The fiber of $\nu$ over $(\C{L},\C{E})$ classifies embeddings of vector bundles 
$\eta:\C{L}\hra\C{E}$, commuting with $F$-structures and 
preserving $D$-level structures.  In particular, this means that the $D$-level structure of $\C{E}$ 
induces an isomorphism between $\C{L}_{|D\times S}$ and the first summand of $\C{O}_{D\times S}^n$. 
(Here we use the assumption that $[g]=1$). 

Consider first an  a priori slightly bigger stack $\C{B}'$ equipped with a morphism \break 
 $\nu':\C{B}'\to FBun^{d}_{GL_1,D,n,\bar{k}}\times_{(X\sm D)^n}\C{A}$, whose fibers classify the same 
data as $\nu$, but $\eta(\C{L})$ is just a subsheaf of $\C{E}$ and not necessary a 
subbundle. As $|D|>l_1-d$, our choice of $l_1$ implies that such an 
$\eta$ is at most unique, therefore $\nu'$ is a closed embedding. 
Since   $FBun^{d}_{GL_1,D,n,\bar{k}}$
is finite and \'etale over $(X\sm D)^n$, this shows that $\C{B}'$ is finite and unramified over
$\C{A}$. 

Next we will show that an open embedding $\C{B}\hra\C{B}'$ is actually an isomorphism.
Assuming the contrary, there exists an $\overline{\fq}$-point of $\C{B}'$ such that 
the corresponding $\C{L}$ is not a subbundle of $\C{E}$ 
(here we identify $\C{L}$ with $\eta(\C{L})$).
Let $\C{L}'$ be the line subbundle of $\C{E}$, containing $\C{L}$. Then $\C{L}'=\C{L}(D')$
for certain finite non-empty subscheme $D'$ of $X\otimes_{\fq}\overline{\fq}$. 
As $\C{E}\in\C{A}(\overline{\fq})$,
we get $|D'|\leq l_1-d$. Also we know that $\phi$ induces an
isomorphism between restrictions of $\C{L}'$ and ${}^{\tau}\C{L}'$
to $ X\otimes_{\fq}\overline{\fq}-\{x_1,\ldots,x_n\}$. Hence there
exist integers $r_1,\ldots,r_n$ such that the divisor ${}^{\tau}(D')-D'$
is equal to $\sum_{i=1}^n r_i x_i$. 

As $\eta$ preserves $D$-level structures, $D'$
is disjoint from $D$. Hence each point of $D'$ has a degree at
least $l_1-d+1$ over $\fq$. Then for each $y\in D'$, all points
$y,{}^{\tau}y,\ldots,{}^{\tau^{l-d}}y$ are distinct. Since
$|D'|\leq l_1-d$, we thus get that ${}^{\tau}(D')\neq D'$. Choose a point
$y\in D'$, which does not lie in ${}^{\tau}(D')$. Let $r$ be the
smallest positive integer such that ${}^{\tau^r}y\notin D'$. Then
$r\leq l_1-d$ and ${}^{\tau^{r-1}}y\in D'$, hence ${}^{\tau^r}y\in
{}^{\tau}(D')\sm D'$. Then the equality
${}^{\tau}(D')-D'=\sum_{i=1}^n r_i x_i$ implies that both $y$ and 
${}^{\tau^r}y$ belong to $\{x_1,\ldots,x_n\}$, contradicting our
assumption.

It remains to show that the projection $\C{B}\to\C{A}$ is 
injective. If not, then there exists a geometric point of $\C{A}$ and
two different degree $d$ line subbundles $\C{L}_1$ and $\C{L}_2$ of the corresponding
vector bundle $\C{E}$, whose restrictions to $D$ coincide. Let
$\wt{\C{E}}$ be the rank two subbundle of $\C{E}$, generated by
$\C{L}_1$ and $\C{L}_2$. Then we get a contradiction $\deg \wt{\C{E}}\geq
\deg\C{L}_1+\deg\C{L}_2+|D|=2d+|D|>l_2$. 
\end{proof}

\section{Cuspidal part of the cohomology}

\begin{Not} \label{N:pure}
a) Let $L/\fq$ be a finitely generated field extension. We will
say that a $Gal(\ov{L}/L)$-module $\C{F}/\overline{\B{Q}_l}$ is
{\em mixed of weight $\leq i$} (resp. {\em pure of weight $i$})
if there is a scheme of finite type $Z/\fq$ with
$\fq(Z)=L$ such that $\C{F}$ is the restriction of a mixed of weight $\leq i$ 
(resp. pure of weight $i$) $\overline{\B{Q}_l}$-sheaf on $Z$. In particular,
if $\C{F}$ is a mixed module of weight $\leq i$, it has a finite weight filtration
$0=\C{F}_j\subset\C{F}_{j+1}\subset\ldots\subset\C{F}_i=\C{F}$ such that each graded piece
$gr_k(\C{F}):=\C{F}_k/\C{F}_{k-1}$ is pure of weight $i$.

b) Let $Y/L$ be a scheme of finite type, and let
$\C{F}$ be an object of $D^b(Y,\overline{\B{Q}_l})$. We will say
that $\C{F}$ is a {\em mixed complex of weight $\leq 0$} if there is a
morphism $\pi:\wt{Y}\to Z$ of schemes of finite type over $\fq$
such that $\fq(Z)=L$, $Y$ is the generic fiber of $\pi$, and
$\C{F}$ is the restriction of a mixed complex on $\wt{Y}$ of weight
$\leq 0$.

c) Combining Deligne's theorem (\cite[6.2.3]{De}) with proper base
change theorem, we obtain that if $\C{F}$ is a mixed complex on $Y$ of
weight $\leq 0$, then each $H^i_c(Y,\C{F})$ is a mixed $Gal(\ov{L}/L)$-module of
weight $\leq i$. We denote by $H^i_{c,\pure}(Y,\C{F})$ the corresponding graded piece $gr_i(H^i_c(Y,\C{F}))$
and will call it pure cohomology with compact support.
\end{Not}

\begin{Rem} \label{R:pure}
Let $Y$ be as in \rn{pure}, and let $\C{F}$ be a mixed complex on $Y$
of weight $\leq 0$.

a) If $f:Y'\to Y$ is a morphism of schemes of finite type over $L$,
then $f^*(\C{F})$ is mixed of weight $\leq 0$. Therefore pure
cohomology $H^{i}_{c,\pure}(Y',f^*(\C{F}))$ is defined.

b) Every open embedding $U\hra Y$ induces an embedding 
$H^i_{c,\pure}(U,\C{F}_{|U})\hra H_{c,\pure}^i(Y,\C{F})$. Indeed,
the kernel of the canonical map $H^i_{c}(U,\C{F}_{|U})\to H_{c}^i(Y,\C{F})$
is a quotient of  $H^{i-1}_{c}(Y\sm U,\C{F}_{|Y\sm U})$, hence by
Deligne's theorem, it is mixed of weight $\leq (i-1)$.

c) Every proper morphism $f:Y'\to Y$ induces a morphism
\[
f^*_{\pure}:H_{c,\pure}^i(Y,\C{F})\to H^i_{c,\pure}(Y',f^*(\C{F})).\]
In particular, we have a restriction map $H_{c,\pure}^i(Y,\C{F})\to
H^i_{c,\pure}(Y',\C{F}_{|Y'})$ for each closed subscheme $Y'$ of $Y$.
Moreover, an argument, similar to b) shows that if $Y_1,\ldots,
Y_k$ are the set of all irreducible components of $Y$, then the natural
restriction map $H_{c,\pure}^i(Y,\C{F})\to\bigoplus_{j=1}^k
H^i_{c,\pure}(Y_j,\C{F}_{|Y_j})$ is an embedding.

d) If $Y/L$ is proper, then $H^i_c(Y,\IC_Y)=H^i(Y,\IC_Y)$ is pure of weight $i$, and thus \break
$H^i_{c,\pure}(Y,\IC_{Y})=H^i_{c}(Y,\IC_{Y})$. More generally, if
$Y$ is an open subscheme of a proper scheme $Y'/L$, then by b),
$H^i_{c,\pure}(Y,\IC_{Y})$ is the image of the natural morphism \break
$H_{c}^i(Y,\IC_{Y})\to H^i(Y',\IC_{Y'})$. 

\end{Rem}
%
%

\begin{Not} \label{N:dm}
Let a Deligne--Mumford stack $Y$ be the quotient of a quasi-projective scheme $X$ by
a finite group $G$, let $Z$ be the quotient $X$ by $G$ in the
category of schemes, and let $q:Y\to Z$ be the natural map.
Following \cite[App. A]{La} we define
$H^i_c(Y,\C{F})$ to be $H^i_c(Z,q_*(\C{F}))$ for every $\C{F}\in
D^b(Y,\overline{\B{Q}_l})$.
\end{Not}

\begin{Rem} \label{R:equiv}
In the notation of \ref{N:dm}, we have

a) $q_*(\IC_Y)=\IC_Z$ (compare \cite[proof of Prop. A.5]{La} and c)
below), therefore $H^i_c(Y,\IC_Y)=H^i_c(Z,\IC_Z)$;

b) If $\C{F}$ is a mixed complex on $Y$ of weight $\leq 0$, then
$q_*(\C{F})$ is mixed of weight $\leq 0$ as well, so we can define
$H^i_{c,\pure}(Y,\C{F}):=H^i_{c,\pure}(Z,q_*(\C{F}))$ which
automatically satisfies all the properties of \rr{pure};

c) $H^i_c(Y,\IC_Y)=H_c^i(X,\IC_X)^G$. Indeed, let $\pi:X\to Z$ be the
quotient map. Then $\pi$ is finite, thus $\pi_*(\IC_X)$ is a
perverse sheaf, and $\IC_Z=\pi_*(\IC_X)^G$. Since
$H_c^i(X,\IC_X)=H_c^i(Z,\pi_*(\IC_X))$, the statement follows from
the fact that any additive functor (e.g., $H^i_c(Z,\cdot)$) between
$\B{Q}$-linear abelian categories commutes with taking invariants
with respect to a finite group.
\end{Rem}

\begin{Emp} \label{E:cusp}
\begin{proof}[Proof of \rp{cusp}]
a) 
First we will show that the space $H^i_{\cusp}$ is
contained in $H_{\pure}^{i,d(\ov{\om})\rho}$. Let $h$ be any element of
$H_{\pure}^i\sm H_{\pure}^{i,d(\ov{\om})\rho}$. Choose sufficiently
large $\mu$ and $D$ such that $h$
belongs to $H_{\pure,D}^{i,\mu}\sm H_{\pure,D}^{i,d(\ov{\om})\rho}$. Hence
$h$ does not vanish on $J\bs
[FBun_{D,n,\overline{\om}}^{\leq\mu}\sm FBun_{D,n,\overline{\om}}^{\leq d(\ov{\om})\rho}]$. 
By \rr{pure} c), $h$ therefore does not vanish on a certain 
irreducible component $C^0$ of 
$J\bs[FBun_{D,n,\overline{\om}}^{\leq\mu}\sm FBun_{D,n,\overline{\om}}^{\leq d(\ov{\om})\rho}]$. 
By \rr{pure} b), $h$ then does not vanish on the closure 
$C$ of $C^0$ in $J\bs FBun_{D,n,\overline{\om}}$. As $C$ is an orispheric substack (use
\rt{red} and \rp{par1}), we see that $h\notin H^i_{\cusp}$. This shows that
$H^i_{\cusp}\subset H_{\pure}^{i,d(\ov{\om})\rho}$.

As the set of orispheric substacks  
is $\ga$-invariant,  $H^i_{\cusp}$ is a $\ga$-invariant subspace of
$H_{\pure}^i$. Therefore  $H^i_{\cusp}$ is contained in $\cap_{g\in\ga}
g(H_{\pure}^{i,d(\ov{\om})\rho})$, hence in $\cap_{g\in\ga}
g(H_{\pure}^{i,\mu})$ for each $\mu\geq d(\ov{\om})\rho$.

Conversely, assume that some $h\in H_{\pure}^i$ does not lie in
$H^i_{\cusp}$. Then the restriction of $h$ to some orispheric
substack $C\subset J\bs FBun_{*,n,\overline{\om}}$ is non-trivial.
Let ${s}$ be the geometric generic point of $C$. It remains to show that for every $\mu$
there exists $g\in\ga$ such that $g({s})\notin J\bs
FBun_{*,n,\overline{\om}}^{\leq\mu}$. Indeed, this would imply that
$g(C)\cap J\bs FBun_{*,n,\overline{\om}}^{\leq\mu}=\emptyset$,
hence $g^{-1}(h)\notin H_{\pure}^{i,\mu}$, thus $h\notin
g(H_{\pure}^{i,\mu})$, as claimed.

By \rp{par1}, ${s}$ can be considered as a geometric point of
$FBun_{P,*,n}$ for certain  maximal parabolic $P$ of
$G$. Moreover, replacing ${s}$ by its $G(\B{O})$-translate, we
may assume that it lies in $FBun^{d,\bar{k},[g]}_{P,*,n}$ with $g=1$.
Let $\la$ be the quasi-fundamental weight of 
$G$ corresponding to $P$.
If for each $g\in P(\B{A})$ we denote by $\C{P}_g$ the fiber at $g(s)$ of the universal $P$-bundle, 
then $\deg(\C{P}_g)_{\la}=\deg(\C{P}_1)_{\la}+\log_q|\la(g)|$. In
particular, $g({s})$ does not belong to 
$J\bs FBun_{*,n,\overline{\om}}^{\leq\mu}$ if $|\la(g)|$ is sufficiently
large.

b) Let $U\subset G(\B{O})\subset G(\B{A})$ be any compact open subgroup. 
By a), the space of invariants $(H^i_{\cusp})^U$ is contained in
$(H_{\pure}^{i,d(\ov{\om})\rho})^U$. By \rr{equiv} c), the latter space
is contained in some $H^{i,d(\ov{\om})\rho}_{\pure,D}$. Since this space is clearly 
finite dimensional, the statement follows.
\end{proof}
\end{Emp}

\section{Drinfeld--Lafforgue case} \label{S:laf}

\begin{proof}[Proof of \rt{laf}]
 
During the proof we will use Drinfeld's notation $FSh_r$ instead of the generic fiber  of $FBun_{GL_r}$. 
Since $FSh_r$ is smooth of relative dimension $2(r-1)$ over $F^{(2)}$,
we have  $H_c^0(FSh_r,\IC)=H_c^{2(r-1)}(FSh_r,\overline{\B{Q}_l}(r-1))$.

\begin{Not} \label{N:rnegl}
Following Lafforgue (\cite{La}), we will say that a representation $V$
of $\Gm_{F^{(2)}}$ is {\em $r$-negligible} if each irreducible
subquotient $V'$ of $V$ is isomorphic to a subquotient of  
the exterior product $V'_1\pp V'_2$ of $(\Gm_F)^2$ 
(composed with the projection $\Gm_{F^{(2)}}\to(\Gm_F)^2$) for some
representations $V'_1$ and $V'_2$ of $\Gm_F$ of dimensions
strictly less than $r$.
\end{Not}


\begin{Lem} \label{L:negl}

For  each $d(\ov{\om})\rho\leq\mu\leq\mu'$ and each
 $D$, the quotient $H_{\pure,D}^{0,\mu'}/H_{\pure,D}^{0,\mu}$ is $r$-negligible.
\end{Lem}

\begin{Rem} \label{R:negl} 
a) Lafforgue proved this assertion under the assumption
that $\mu$ is sufficiently large as a function of $D$ (see
\cite[Cor. VI. 21]{La}).

b) Though Lafforgue \cite{La} and Drinfeld \cite{Dr} worked only in the case when
$J\subset\B{A}\m$ is a cyclic group generated by an element of
degree one, the general case applies without any changes.
\end{Rem}

First we will show the theorem, assuming the lemma.

a) \rl{negl} implies that   
$H_{\pure}^0/H_{\pure}^{0,d(\ov{\om})\rho}$ is $r$-negligible. By \rp{cusp} a), 
$H_{\pure}^0/H^0_{\cusp}$ is a subrepresentation of
$\bigoplus_{g\in G(\B{A})}H_{\pure}^0/g(H_{\pure}^{0,d(\ov{\om})\rho})$. Since
the latter space is isomorphic as a $\Gm_{F^{(2)}}$-representation to the direct sum of
$G(\B{A})$ copies of $H_{\pure}^0/H_{\pure}^{0,d(\ov{\om})\rho}$, we get that
$H_{\pure}^0/H^{0}_{\cusp}$ is $r$-negligible as well.

By (one of) the Main result(s) of Lafforgue (see \cite[Lem. IV. 25
and Thm. VI. 27]{La}), there exists an exhausting $\ga\times\Gm_{F^{(2)}}$-invariant 
filtration $0=\wt{V}_0\subset\wt{V}_1\subset\wt{V}_2\subset\ldots\subset
H^0$ such that each quotient $\wt{V}_{2i+1}/\wt{V}_{2i}$ is  
$r$-negligible, and $\bigoplus_{j}\wt{V}_{2j}/\wt{V}_{2j-1}$ 
is isomorphic to $\bigoplus_{\pi} (\pi\pp\rho_{\pi}\pp\check{\rho}_{\pi})$ (compare
\rr{negl} b)).
Actually Lafforgue has shown that the last isomorphism holds after semi-simplification.
However, it was observed by Drinfeld that there are no non-trivial extensions 
between non-isomorphic cuspidal representations $\pi_1$ and $\pi_2$ of $GL_r(\B{A})$,
therefore the former representation is automatically semi-simple.
(Indeed, choose a compact and open subgroup $U$ of $GL_r(\B{A})$ such that
 $(\pi_1)^U\neq 0$ and $(\pi_2)^U\neq 0$. It remains to show that there are no
 non-trivial extensions between $(\pi_1)^U$ and $(\pi_2)^U$ as  representations of the Hecke algebra 
$\C{H}_U:=\C{H}(GL_r(\B{A}),U)$. Since the center of $\C{H}_U$ contains $\C{H}(GL_r(F_v),GL_r(\C{O}_v))$ for 
almost all points $v$ of $X$, the strong multiplicity one theorem for $GL_r$ implies that
 $(\pi_1)^U$ and $(\pi_2)^U$ has non-equal central characters of $\C{H}_U$. This implies the assertion).

Since we know that 
$\rho_{\pi}\pp\check{\rho}_{\pi}$ is pure (Ramanujan
conjecture \cite[Thm. VI.10]{La}) and that $H^0_{\pure}/H^0_{\cusp}$
is $r$-negligible, the induced filtration of $H^0_{\cusp}$ satisfies the required property.
(By an induced filtration we mean the filtration $\{V_j\}_j$ such that 
each $V_j$ is the intersection of  $H^0_{\cusp}$ with the image of $\wt{V}_j$ in $H^0_{\pure}$.)

b) In the case $r=2$, Drinfeld (\cite{Dr}) constructed a
compactification $\overline{J\bs FSh^{\leq\mu}_{2,D}}$ of each
$J\bs FSh^{\leq\mu}_{2,D}$ with $\mu$  sufficiently large (which is the
coarse  moduli space of the compactification considered by
Lafforgue \cite[III]{La}). Moreover, for each $\mu_1\leq\mu_2$ and
$D_1\subset D_2$, there exists a canonical
morphism $\overline{J\bs
FSh^{\leq\mu_2}_{2,D_2}}\to\overline{J\bs
FSh^{\leq\mu_1}_{2,D_1}}$, and the inverse limit $\overline{J\bs FSh_2}$
of the $\overline{J\bs FSh^{\leq\mu}_D}$'s is equipped with an
action of the adelic group $G(\B{A})$.

Though compactifications $\overline{J\bs FSh^{\leq\mu}_{2,D}}$
are not smooth, their middle cohomology $H^2(\overline{J\bs
FSh^{\leq\mu}_{2,D}},\overline{\B{Q}_l}(1))$ is pure of weight zero and satisfies Poincar\'e
duality. As a result, their direct limit
$H^2(\overline{J\bs FSh_2},\overline{\B{Q}_l}(1))$ is equipped with a
non-degenerate pairing.
Furthermore, Drinfeld's main result says that the orthogonal
complement $H_{Drin}\subset
H^2(\ov{J\bs FSh_2},\overline{\B{Q}_l}(1))$ of the set of all
Chern classes of orispheric curves is isomorphic to the direct sum
$\bigoplus_{\pi}(\pi\pp\rho_{\pi}\pp\check{\rho}_{\pi})$. By the Poincar\'e 
duality, $H_{Drin}$ can be described as a set of all elements of
$H^2(\overline{J\bs FSh_2},\overline{\B{Q}_l}(1))$, vanishing on all
orispheric curves.

Since $H^2(\overline{J\bs FSh_2},\overline{\B{Q}_l}(1))$ is pure of weight zero, we get 
a canonical embedding $H^0_{\cusp}\hra H^2(\overline{J\bs FSh_2},\overline{\B{Q}_l}(1))$, which 
therefore induces an embedding $H^0_{\cusp}\hra H_{Drin}$. Combining this with the result of a),
we thus conclude that $H^0_{\cusp}=H_{Drin}$.
\end{proof}

\begin{Emp}  \label{E:negl}
\begin{proof}[Proof of \rl{negl}]
 By \rr{pure}, c), it suffices to check that
$H_{c,\pure}^{2(r-1)}(S,\overline{\B{Q}_l})$ is $r$-negligible for
every irreducible component $S$ of $J\bs [FSh_{r, D}^{\leq\mu'}\sm FSh^{\leq\mu}_{r,D}]$.
Our strategy will be to reduce the statement to the case,
where $S$ is an $\B{A}^{r-m}$-bundle over an open substack $FSh'_{m,D}$ of
$FSh_{m,D}$ for certain $m<r$. Indeed, in this case
$H_{c,\pure}^{2(r-1)}(S, \overline{\B{Q}_l})\cong H_{c,\pure}^{2(m-1)}(FSh'_{m,D},\overline{\B{Q}_l})$ 
would be isomorphic to a subspace of $H_{c,\pure}^{2(m-1)}(FSh_{m,D},\overline{\B{Q}_l})$
 (use \rr{pure} b)). Hence by \cite[Prop. VI.15]{La} it is $(m+1)$-negligible,
hence $r$-negligible.

By \rt{red} and \rp{par1}, $S$ is an open substack of an
orispheric substack. Let $P$ be the  parabolic subgroup of
$G$, corresponding to $S$. Since the statement for large $D$'s
implies that for small ones, we may increase $D$ during the proof.
Hence using  \rp{par1} and \rr{pure} b), we may replace $S$ by the
generic fiber of $FSh_{r,P,D}$.  In other words it will suffice
to show the statement for (each connected component of the) stack
$S$ classifying pairs consisting of an $F$-sheaf of rank $r$ with
$D$-level structure $(\C{E},\psi;x_1,x_2;\phi)$ and a subbundle
$\C{A}$ of $\C{E}$ of rank $m$ (determined by $P$) such that
$\phi(\C{A})\subset{}^{\tau}\C{A}(x_1)$ and  $D$-level structure $\psi$ maps $\C{A}_{|D}$
into the first $m$-coordinates.

Applying if necessary transformation sending an
$F$-sheaf to its dual, we may assume that $S$ classifies pairs for
which $\C{B}:=\C{E}/\C{A}$ is a pullback of a vector bundle over
$X$, and $\phi$ induces a canonical isomorphism
$\C{B}\isom{}^{\tau}\C{B}$. Thus we can replace $S$ with its open and 
closed substack, over which both the quotient $\C{B}=\C{E}/\C{A}$ and its induced 
$D$-level structure are constant. Moreover, we may fix a coweight $\nu$ of $SL_m$ 
and an integer $d$ and to replace $S$ by its open substack such that 
$\C{A}$ belongs to $FSh^{\leq\nu;d}_{m,D}$. 
Enlarging $|D|$, we may assume that $\Hom(\C{B},\C{A}'(-D))=0$ for each
$({}^{\tau}\C{A}\hra \C{A}'\hookleftarrow\C{A})\in
FSh^{\leq\nu;d}_{m}$.

Choose a  complete flag of subbundles
$0=\C{B}_0\subset\C{B}_1\subset\ldots\subset\C{B}_{k-m-1}\subset\C{B}_{k-m}=\C{B}$
defined over $\fq$, and put $M:=d-2g+1-|D|-\max_j
\deg(\C{B}_j/\C{B}_{j-1})$. By induction, we may assume the
statement for smaller $m$'s (the assumption is vacuously true for
$m=1$), and thus to replace $S$ by an open substack $S'$
characterized by the condition that $\C{A}$ does not have
$F$-subsheaves of rank $(m-1)$ and degree larger than $M$. We claim that
the natural forgetful map $\pi:S'\to FSh^{\leq\nu;d}_{m,D}$ is smooth, and all its non-empty
fibers are isomorphic to $\B{A}^{r-m}$. 

As the smoothness statement was shown in \cite[II, Thm. 11]{La1}
(and easily follows from our arguments below), it remains to
show the statement about fibers. Denote the image of $\pi$ by
$FSh'_{m,D}$, and let $s$ be any  geometric point of $FSh'_{m,D}$.
We claim that the fiber $\pi^{-1}(s)$ is canonically isomorphic to
the space of extensions $\wt{\C{E}}\in\Ext(\C{B},\C{A}_s(-D))$
whose image in $\Ext(\C{B},\C{A}'_s(-D))$, induced by the embedding $\C{A}_s\hra\C{A}'_s$,
coincides with the image of  ${}^{\tau}\wt{\C{E}}\in\Ext(\C{B},{}^{\tau}\C{A}_s(-D))$, 
induced by the embedding ${}^{\tau}\C{A}_s\hra\C{A}'_s$. Indeed, given
$(\C{E}_s,\psi)$ in $\pi^{-1}(s)$ we can get $\wt{\C{E}}$ as the
kernel of the composition of the projection $\C{E}_s\to
\C{E}_{s|D}$, the level structure $\psi$, and the projection to
the first $m$ factors. Conversely, given $\wt{\C{E}}$, we can
define $\C{E}_s$ as $(\wt{\C{E}}\oplus\C{A}_s)/\C{A}_s(-D)$, and
the first $m$ (resp. the last $r-m$) coordinates of $\psi$ to be the
composition of the projection
$\C{E}_s\to\C{A}_s/\C{A}_s(-D)=\C{A}_{s|D}$ (resp.
$\C{E}_s\to\wt{\C{E}}_{|D}\to\C{B}_{|D}$) and $D$-level structure
of $\C{A}$ (resp. $\C{B}$).

Put $U:=\Ext(\C{B},\C{A}_s(-D))$ and $V:=\Ext(\C{B},\C{A}'_s(-D))$.
Then embeddings $\C{A}_s\hra\C{A}'_s$ and
${}^{\tau}\C{A}_s\hra\C{A}'_s$ define a linear  homomorphism
$\la:U\to V$ and a $\tau$-linear homomorphism $\psi:U\to V$, respectively. 
As $\C{A}'/\C{A}$ is supported at one point, $\la$ is surjective, and
our assumption $\Hom(\C{B},\C{A}'_s(-D))=0$ implies that its kernel
is of dimension $\rk(\C{B})=k-m$. As $\pi^{-1}(s)$ is isomorphic to
$\Ker(\la-\psi)$, \cite[II, Lem. 18 iii)]{La1} and its proof
(compare \cite[\S 4, Lem. 1]{Dr1}) implies that in order to show the
statement, it remains to check that there is no linear functional
$L$ on $V$ such that $L(\psi(u))=L(\la(u))^q$ for each $u\in U$.

Assume that such an $L$ exist. By Serre duality, $\check{V}\cong
H^0(X,K_X\otimes\C{B}\otimes\check{\C{A}}'(D))$, so $L$ defines
a non-trivial $\phi$-equivariant section of
$K_X\otimes\C{B}\otimes\check{\C{A}}(D)$. Hence the exists $j$
such that $L$ defines a non-trivial $\phi$-equivariant section of
$K_X\otimes(\C{B}_j/\C{B}_{j-1})\otimes\check{\C{A}}(D)$, thus a
non-trivial $\phi$-equivariant morphism $\C{A}\to
K_X\otimes(\C{B}_j/\C{B}_{j-1})(D)$. Its kernel is a $F$-subsheaf
of $\C{A}'$ rank $l-1$ and degree at least
$d-(2g-2)-\deg(\C{B}_j/\C{B}_{j-1})-|D|>M$, contradicting our
definition of $S'$. This completes the proof of the claim, which (as it was observed in the 
beginning of the proof) implies the lemma. 
\end{proof}
\end{Emp}
\appendix
\section{}

In the appendix we give proofs (or sketches of proofs) of basic properties
of $G$-bundles, affine grassmannians and the stacks of Hecke, used in the paper. 
Though at least some of these facts are considered well known among experts, 
we were not able to find any reference in the literature, so
we sketch the argument for completeness.

\begin{Emp} \label{E:conn}
\begin{proof}[Proof of \rl{conn}]
In the case $G=\B{G}_m$, the required isomorphism between
$\pi_0(Bun_{\B{G}_m})=\pi_0(Pic)$ and $X_*(\B{G}_m)=\B{Z}$ is given by the degree map.
This implies the statement in the case of torus. 
In general, the embedding $T\hra G$ defines a surjective map
$\pi: X_*(T)=\pi_0(Bun_T)\to\pi_0(Bun_G)$ (see \cite[App.]{DS}).
Moreover, using standard reduction 
to the $GL_2$-case (see \cite[App.]{DS}, where the semi-simple simply-connected 
case is treated, and compare the proof of \rp{affgr} b)), we get that $\pi$ factors through
$\pi_1(G)$. Thus it remains to show the injectivity of the resulting canonical map
$\ov{\pi}_G:\pi_1(G)\to\pi_0(Bun_G)$.

When  $G^{\der}$ is simply connected, the injectivity of $\ov{\pi}_G$ follows from the fact
that the composition map $\pi_1(G)\to\pi_0(Bun_G)\to\pi_0(Bun_{G^{\ab}})=\pi_1(G^{\ab})$
is an isomorphism. For a general $G$, choose a central extension $\nu:H\to G$
such that $H^{\der}$ is simply-connected, and $S:=Ker(\nu)$ is a
split torus (see \cite[Prop. 3.1]{MS}). Then we have the following 
commutative diagram, induced by $\nu$:

$$
\CD
\pi_1(H)              @>\ov{\pi}_H >> \pi_0(Bun_H)\\
@V\pi_1(\nu)VV                  @V{\pi_0(\nu)}VV\\
\pi_1(G)              @>\ov{\pi}_G >> \pi_0(Bun_G)
\endCD
$$
As $Bun_S$ acts transitively on all geometric fibres of the projection
$Bun_H\to Bun_G$, 
$\pi_0(Bun_S)=\pi_1(S)$ acts transitively on all fibers of
$\pi_0(\nu)$. Since  $\ov{\pi}_H$ is injective and $\pi_1(S)$-equivariant, and since
$\pi_1(G)=\pi_1(S)\bs\pi_1(H)$, 
the injectivity of $\ov{\pi}_G$ follows.
\end{proof}
\end{Emp}
\begin{Lem} \label{L:par}
Let $P$ be a  parabolic subgroup of $G$, and let $\C{G}$
be a $G$-bundle over a scheme $S$. Then there exists a canonical
bijection between:

(i) $P$-structures of $\C{G}$;

(ii) families of line subbundles $\{\C{A}_{\la}\subset
\C{G}_{\la}\}_{\la}$ (indexed by  characters of $P$ which are
dominant weights of $G$) such that
$\C{A}_{\la_1+\la_2}=\C{A}_{\la_1}\otimes\C{A}_{\la_2}$ 
(considered as subbundles of 
$\C{G}_{\la_1+\la_2}\subset\C{G}_{\la_1}\otimes\C{G}_{\la_2}$) for
each $\la_1$ and $\la_2$;

(iii) families of line subbundles $\{\C{A}_{\la_i}\subset
\C{G}_{\la_i}\}_{i}$ (indexed by characters of $P$ which are
quasi-fundamental weights of $G^{\ad}$) satisfying the Pl\"ucker
relations: for each tuple of non-negative integers $\{k_i\}_i$, the line subbundle
$\otimes_{i}(\C{A}_{\la_i})^{k_i}\subset\otimes_{i}(\C{G}_{\la_i})^{\otimes
k_i}$ is contained in $\C{G}_{\sum_i k_i\la_i}$.
\end{Lem}
\begin{proof}

(i)$\implies$(ii). Let $\C{P}$ be a $P$-structure of $\C{G}$. If a dominant weight $\la$ of $G$ 
is a character of $P$, then the highest weight line $l_{\la}$ of $V_{\la}$ is a $P$-subrepresentation 
of $V_{\la}$. Moreover, line subbundles $\C{A}_{\la}:=P\bs[\C{P}\times l_{\la}]\subset
P\bs[\C{P}\times V_{\la}]=\C{G}_{\la}$ satisfy (ii).

(ii)$\implies$(iii) is clear.

(iii)$\implies$(i). Suppose that we are given a family
$\{\C{A}_{\la_i}\subset \C{G}_{\la_i}\}_{i}$, satisfying the
Pl\"ucker relations. We claim that this family defines a canonical
$P$-structure of $\C{G}$ or what is the same a canonical section
of the fibration $P\bs\C{G}\to S$. By the uniqueness assertion, the
statement is local in the \'etale topology on $S$, so we may assume
that $\C{G}$ is trivial. Hence $P\bs\C{G}\cong (P\bs G)\times S$,
thus we want to construct an $S$-point of $P\bs G$. We have a
canonical embedding $P\bs G\hra\prod_i \B{P}(V_{\la_i})$.
Each $\C{A}_{\la_i}\subset \C{G}_{\la_i}$ corresponds to an
$S$-point of $\B{P}(V_{\la_i})$, and the fact that the
$\{\C{A}_{\la_i}\subset \C{G}_{\la_i}\}_{i}$ satisfy the
Pl\"ucker relations means precisely that the corresponding
$S$-point of the product $\prod_i \B{P}(V_{\la_i})$ defines an $S$-point 
of $P\bs G$.
\end{proof}


\begin{Lem} \label{L:open}
$Bun_G^{\leq\mu}$ is an open substack of $Bun_G$.
\end{Lem}
\begin{proof}
As any dominant weight of $G^{\ad}$ is a linear combination of
quasi-fundamental weights with rational non-negative coefficients,
it is enough to check the condition $\deg (\C{B}_s)_{\la}\leq\langle \mu,\la\rangle $ 
only when $\la$ is a quasi-fundamental weight. Since the number of quasi-fundamental weights is finite,
it remains to check the openness of this condition for a given 
quasi-fundamental weight $\la_i$. 

Let $P_i\supset B$ be the maximal parabolic subgroup corresponding to $\la_i$; then 
$(\C{B}_s)_{\la_i}$ depends only on $P_i\times_B\C{B}_s$. By \rl{par}, 
condition $\deg (\C{B}_s)_{\la_i}\leq\langle \mu,\la_i\rangle $ is equivalent to the assertion
that $(\C{G}_{{s}})_{\la_i}$ has no line subsheaves of degree $\langle \mu,\la_i\rangle +1$, satisfying Pl\"ucker
relations. Since Pl\"ucker relations are closed, and since line subsheaves of given degree 
of a given vector bundle are represented by a projective scheme, our condition  
$\deg (\C{B}_s)_{\la_i}\leq\langle \mu,\la_i\rangle $ is therefore open, as claimed. 
\end{proof}

\begin{Emp} \label{E:repr}
\begin{proof}[Proof of \rl{repr}]

a) First we consider the case of $G=GL_n$. Fix an ample line bundle $\C{O}(1)$
an $X$. Then given $\mu$, there exists $m\in\B{N}$ such that for every
$S/\fq$ and every $\C{E}\in Bun^{\leq\mu}(S)$, we have:

i) the direct image $\pr_*(\C{E}(m))$ is a vector bundle on $S$,

ii) $R^1\pr_*(\C{E}(m))=0$, and 

iii) $\C{E}$ is a quotient of $\pr^*[\pr_*(\C{E}(m))](-m)$. 

\noindent Next given ($\mu$ and)
$m$, we get that $\pr_*(\C{E}(m))$ is a subbundle of
$pr_*(\C{E}(m)_{|D\times S})\isom pr_*(\C{O}^n_{|D\times S})$ for
each $|D|$ sufficiently large, where the last isomorphism is induced by 
the $D$-level structure. Thus the functor $(\C{E},\psi)\mapsto (\pr_*(\C{E}(m)),\C{E})$
embeds $Bun_{D}^{\leq\mu;\nu}$ into a quasi-projective scheme 
classifying pairs consisting of a subbundle $\C{H}$ of $pr_*(\C{O}^n_{|D\times S})$ of given rank,
and a locally free quotient of $[pr^*\C{H}](-m)$ of rank $n$ and
degree $\nu$. This implies the representability of our functor in
the $GL_n$-case.

For the general case, we can proceed as in (\cite[Sec.4]{Be1}): 
Choose an embedding of $G$ into $GL_n$. Since $Bun_{G}^{\leq \mu;\nu}$ 
is quasi-compact, the map $\C{G}\mapsto GL_n\times_G\C{G}$
maps $Bun_{G,D}^{\leq\mu;\nu}$ into some
$Bun_{GL_n,D}^{\leq\mu';\nu'}$. Next, using the fact that
$G$ is reductive and, therefore, $GL_n/G$ is affine, we conclude that
the induced map $Bun_{G,D}^{\leq\mu;\nu}\to Bun_{GL_n,D}^{\leq\mu';\nu'}$ is affine, 
implying the representability assertion. Finally, since coherent sheaves on curves have no second 
cohomology, the smoothness assertion follows from deformation theory.

b) First, we claim that there exists a faithful representation $V$ of $G$, which is a direct
sum of Weyl modules. Indeed, our claim is equivalent to the assertion that
the intersection of kernels of all Weyl modules is trivial. As $G$ is reductive, this 
intersection is obviously central, so it is contained in the maximal torus $T$
of $G$. But dominant weights of $G$ generate the group of all characters of $T$, so their kernels 
have a trivial intersection, as claimed.

Now our strategy will be very similar to that of \cite[A.5]{Ga}, where the
corresponding statement for $Gr_{\om}$, defined below, is shown. Choose  $N\in\B{N}$ such that
$-N\leq\langle \om_i,\la\rangle \leq N$ for each $i=1,\ldots,n$ and each weight $\la$ of $V$.
Consider a substack $Hecke'$ of $Hecke_{n}$, consisting of tuples
$(\C{G},\C{G}'; x_1,\ldots,x_n; \phi)\in Hecke_{n}$ satisfying

i) $\C{G}'_V(-N(\sum_{i=1}^n
\Gm_{x_i}))\subset\phi(\C{G}_V)\subset\C{G}'_V(N(\sum_{i=1}^n
\Gm_{x_i}))$ and 

ii) $\det(\C{G}_V)=\det(\C{G}'_V)(\sum_{i=1}^n
\langle \det(V),\om_i\rangle \Gm_{x_i}).$\\
Then $Hecke_{n, \overline{\om}}$ is a closed substack of
$Hecke'$ (use \rr{inv} c)), so it remains to show the representability and
projectivity of $Hecke'\to Bun\times X^n$. 

Consider another stack $Hecke''$ classifying data 
$(\C{E},\C{G}'; x_1,\ldots,x_n; \phi)$, where \break 
$(\C{G}'; x_1,\ldots,x_n)\in (Bun\times X^n)(S)$, 
$\C{E}$ is a vector bundle over $X\times S$ of rank $\dim V$, and $\phi$ is an isomorphism between the restrictions
of $\C{E}$ and $\C{G}'_V$ to $(X\times S)\sm(\Gm_{x_1}\cup\ldots\cup\Gm_{x_n})$ satisfying 
conditions i) and ii) as above with $\C{G}_V$ replaced by $\C{E}$.
Then $Hecke''$ is represented by a closed substack of a relative grassmannian over
$Bun\times X^n$, hence  $Hecke''\to Bun\times X^n$ is representable and projective.

Let $\wt{\C{G}}$ be the universal $GL(V)$-bundle over
$Hecke''\times X$, corresponding to $\C{E}$. Then the quotient
$G\bs\wt{\C{G}}$ has a canonical section $l$ over $(Hecke''\times
X)\sm (\Gm_{x_1}\cup \ldots \cup \Gm_{x_n})$, corresponding to the universal 
isomorphism of $GL(V)\times_G\C{G}'$ with $\wt{\C{G}}$ over this
set. Moreover, $Hecke'$ is the largest substack $\C{A}\subset Hecke''$ such that
$l$ extends to a regular section on all of $\C{A}\times X$. 
We want to show that $Hecke'$ is a closed substack of $Hecke''$.

As $G$ is reductive, $G\bs GL(V)$ is affine. Thus $G\bs\wt{\C{G}}$ is affine over $Hecke''\times X$.
Since the question is local in the Zariski topology on both  $Hecke''$ and $X$, 
we are thus reduced to the following assertion:
Suppose we are given a scheme $S$, $n$ points $x_1,\ldots,x_n\in X(S)$ and a regular 
function $l$ on $U:=(S\times X)\sm(\Gm_{x_1}\cup\ldots\cup\Gm_{x_n})$. Then 
the functor $\C{F}_l/S$, which for every scheme $T$ classifies morphisms $f\in\Hom(T,S)$ such that
$f^*(l)$ extends to a regular function on $T\times X$, is represented by a closed subscheme of $S$.

Let $j:U\hra S\times X$ be the open embedding, let 
$s\in\Gm(X\times S,j_*(\C{O}_{U})/\C{O}_{X\times S})$ be the image of 
$l$, and for each $i=1,\ldots,n$, let $s_i$ be the restriction of $s$ to $\Gm_{x_i}\isom S$.
Then  $f\in\Hom(T,S)$ belongs to  $\C{F}_l(T)$ if and only if $f^*(s_i)=0$ for each $i$. Therefore
$\C{F}_l$ is represented by the intersection of the schemes of zeros of the $s_i$'s.

c) Since 
$Hecke'_{n,\overline{\om}}=Hecke_{1,\om_1}\times_{Bun} Hecke_{1,\om_2}\times_{Bun}\ldots\times_{Bun} 
Hecke_{1,\om_n}$, it is projective over $X^n\times Bun$ (by b)) hence over  $Hecke_{n,\overline{\om}}$. 
The second assertion follows from \rl{red} a) below applied $(n-1)$-times. 
\end{proof}
\end{Emp}

\begin{Not}
Consider a formal disc $\C{D}:=\Spec k[[x]]$ (resp. punctured formal disc $\C{D}^*:=\Spec k((x))$) over $k$.
For each scheme $X/k$, we denote by $X_{\C{D}}$ (resp. $X_{\C{D}^*}$)
be a scheme (resp. functor) over $k$ such that  $X_{\C{D}}(A)=X(A\widehat{\otimes}_k k[[x]])$ 
(resp.  $X_{\C{D}^*}(A)=X(A\widehat{\otimes}_k k((x)))$) for each $k$-algebra $A$. 
In particular,  $X_{\C{D}}$ is just the inverse limit of the   $X_{\C{D}_i}$'s, where  
$\C{D}_i:=\Spec  k[[x]]/(x^i)$, and  $X_{\C{D}_i}=R_{\C{D}_i/k}X$ is the Weil restriction of scalars. 
\end{Not}

\begin{Def} \label{D:locaff}
Let $Gr$ be an ind-scheme, which classifying $G$-bundles 
over a formal disc $\C{D}$, trivialized over a punctured formal disc $\C{D}^*$ (compare \cite{Ga}).
$Gr$ is called a {\em local affine grassmannian}.
As in the global case, for each $\om\in X_*^+(T)$ we define a closed subscheme $Gr_{\om}$ of $Gr$ 
and a locally closed subscheme $Gr^0_{\om}$.
\end{Def}

\begin{Rem}
$G_{\C{D}^*}$ is the group of automorphisms of the trivial $G$-bundle on $\C{D}^*$.
In particular, $G_{\C{D}^*}$ acts naturally on $Gr$, making $Gr$ a homogeneous space. 
Moreover, $Gr_{\om}$ and $Gr^0_{\om}$ are 
$G_{\C{D}}$-invariant subschemes of $Gr$.
\end{Rem}

\begin{Lem} \label{L:red}
a) For each two positive integers $k<n$, put 
$$U_{k,n-k}:=\{(x_1,\ldots,x_n)\in X^{n} | x_i\neq x_j \text{ for each } i\leq k<j\}.$$ 
Then for each $\ov{\om}_1\in X^+_*(T)^k$ and $\ov{\om}_2\in X^+_*(T)^{n-k}$,
the forgetful morphism 
$$\pi: Hecke_{k,\ov{\om_1}}\times_{Bun} Hecke_{n-k,\ov{\om_2}}\to Hecke_{n,(\ov{\om_1},\ov{\om_2})}$$ 
is an isomorphism over $U_{k,n-k}$.
Also over $U_{k,n-k}$, $Hecke_{n,(\ov{\om_1},\ov{\om_2})}$ is canonically isomorphic to the 
fiber product  
\[
Hecke_{(k,\ov{\om_1}),(n-k,\ov{\om_2})}:=(Hecke_{k,\ov{\om_1}}\times Hecke_{n-k,\ov{\om_2}})
\times_{Bun\times Bun} Bun,
\] 
taken with respect to the projection 
$p'\times p':Hecke_{k,\ov{\om_i}}\times Hecke_{n-k,\ov{\om_2}}\to Bun\times Bun$ (see \rn{proj})
and the diagonal morphism $Bun\to Bun\times Bun$. 

b) Over $X^n\sm\Dt$, $Gr_{n,\overline{\om}}$ is canonically isomorphic to the product $\prod_{i=1}^n
Gr_{1,\om_i}$.

c) There exists a smooth surjective morphism $S\to Bun\times X$ with connected fibers 
such that $Hecke_{1,\om}\times_{X\times Bun}S$ is isomorphic over $S$ to $Gr_{1,\om}\times_X S$.

d) Each $Gr_{1,\om}$ is a Zariski locally trivial fibration over
$X$ with fiber $Gr_{\om}$.
\end{Lem}
\begin{proof}
a) Let $(\C{G},\C{G}';y_1,\ldots,y_{n};\phi)$ be an $S$-point of
$Hecke_{n,(\ov{\om_1},\ov{\om_2})}\times_{X^n}U_{k,n-k}$. Define $\C{G}_1$ (resp. $\C{G}_2$)
to be the $G$-bundle over $X\times S$ whose restriction
to the complement of $\cup_{i\leq k}\Gm_{y_i}$ is that of $\C{G}'$ (resp. $\C{G}$), 
restriction to the complement of $\cup_{i>k}\Gm_{y_i}$ is that of $\C{G}$ (resp. $\C{G}'$), 
and the gluing is done with help of $\phi$. 
Then $(\C{G},\C{G}';\ldots)\mapsto (\C{G},\C{G}_2,\C{G}';\ldots)$ gives us the inverse (of the restriction)
of $\pi$, while the map 
\[
(\C{G},\C{G}';\ldots)\mapsto ((\C{G}_1,\C{G}';\ldots),(\C{G}_2,\C{G}';\ldots))
\]
gives us the required isomorphism
\[
Hecke_{n,(\ov{\om_1},\ov{\om_2})}\times_{X^n}{U_{k,n-k}}\to 
Hecke_{(k,\ov{\om_1}),(n-k,\ov{\om_2})}\times_{X^n}{U_{k,n-k}}.
\]

b) follows from the second part of a) applied $(n-1)$-times.

c) The space $S$, classifying pairs $(x,\C{G})\in
X\times Bun$ together with a trivialization of $\C{G}$ over the
completion of the graph of $x$ (or sufficiently large (depending on $\ov{\om}$) 
nilpotent neighborhood of $x$), 
satisfies the required property by \cite{BL}.

d) As explained in \cite[2.1.2]{Ga}, the fibration
$Gr_{1,\om}\to X$ becomes trivial over a certain principal bundle
over $X$, whose structure group $\C{A}$ is the inverse limit of the
$\C{A}_k$'s with $\C{A}_k(R)=Aut_R(R[t]/(t^{k+1}))$. As every
$\C{A}$-bundle is locally trivial in the Zariski topology, the
statement follows.
\end{proof}

\begin{Prop} \label{P:affgr}
a) The reduced schemes $(Gr^0_{n, \overline{\om}})_{\red}$ and  $(Gr'^0_{n, \overline{\om}})_{\red}$ are
smooth over $X^n$ of relative dimension $\langle 2\rho,\sum_{i=1}^n \om_i\rangle $. Moreover,
$Gr^0_{n, \overline{\om}}$ and  $Gr'^0_{n, \overline{\om}}$ are reduced, unless
the characteristic of $k$ is two, and $G$ has a direct factor isomorphic to $PGL_2$ or $PO_{2m+1}$.

b) Both $Gr_{n,\overline{\om}}$ and  $Gr'_{n,\overline{\om}}$ are irreducible.
\end{Prop}
\begin{proof}

a) By \rl{repr} c), 
$Gr^0_{n, \overline{\om}}\cong Gr'^0_{n, \overline{\om}}\times _{X^n} (X^n\sm \Dt)$. Therefore 
it remains to show the statement for $Gr'^0_{n, \overline{\om}}$.
Secondly, as 
$$Gr'^0_{n, \overline{\om}}=Hecke^0_{1,\om_1}\times_{Bun}\ldots\times_{Bun}
Hecke^0_{1,\om_{n-1}}\times_{Bun} Gr^0_{1,\om_n},$$ \rl{locisom}  
and \rl{red} imply that it will suffice to show that 
$(Gr^0_{\om})_{\red}$ is non-singular of dimension $\langle 2\rho,\om\rangle $ and that
$Gr^0_{\om}$ is reduced  unless
$\har k=2$, and $G$ has a direct factor isomorphic to $PGL_2$ or $PO_{2m+1}$.

By the Cartan decomposition, $(Gr^0_{\om})_{\red}$ is a homogeneous
space for the action of $G_{\C{D}}$ of dimension $\langle 2\rho,\om\rangle $. Thus $(Gr^0_{\om})_{\red}$
is non-singular, and the smoothness of $Gr^0_{\om}$ is equivalent to the smoothness at some point.
Assume now either that $\har k\neq 2$ or that 
$G$ does not have a direct factor isomorphic to $PGL_2$ or $PO_{2m+1}$. 

Let $(\C{G},\phi:\C{G}_{|\C{D}^*\times Gr^0_{\om}}\isom
G\times\C{D}^*\times Gr^0_{\om})$ be the universal object over
$Gr^0_{\om}$. Then for each Weyl module $V_{\la}$, $\phi$ gives 
an embedding of $\C{G}_{\la}$ into $triv_{\la}(\langle \la,\omega\rangle x)$, where we write 
$triv_{\la}$ instead of $V_{\la}\times\C{D}\times Gr^0_{\om}$. Let $\wt{\C{B}}$ be the
stack over $Gr^0_{\om}$ classifying $B$-structures of
$\C{G}$, and let us replace the universal object over $Gr^0_{\om}$ by its pullback to $\wt{B}$.
Then each $b\in\wt{\C{B}}$ defines a
line subbundle $\C{L}_{\la}$ of the fiber of $\C{G}_{\la}$ over $X\times\{b\}$. 

Consider a substack $\C{B}'$ of $\wt{\C{B}}$
consisting of those points $b$ such that  for each $\la\in X_*^+(T)$, the corresponding 
line subbundle $\C{L}_{\la}$ is not contained in 
$triv_{\la}((\langle \la,\omega\rangle -1)x)$.
Since the conditions for $\la_1$ and $\la_2$ appearing in the definition of $\C{B}'$ 
imply that for $\la_1+\la_2$, $\C{B}'$ is defined inside $\wt{\C{B}}$ by finitely many open conditions, 
so $\C{B}'$ is open in $\C{B}$. Also since $\wt{\C{B}}$ is pro-smooth over $Gr^0_{\om}$,  
it will suffice to show that $\C{B}'$ is reduced.

For each $b\in\C{B}'$, the corresponding line subsheaf
$\C{L}_{\la}(-\langle \la,\omega\rangle x)$ of $triv_{\la}$ is a
subbundle. Furthermore, these subbundles satisfy the Pl\"ucker
relations, so by \rl{par}, the rule $b\mapsto
\C{L}_{\la}(-\langle \la,\omega\rangle x)$ defines a morphism
$f$ from $\C{B}'$ to the reduced scheme $\C{B}:=(G/B)_{\C{D}}$, classifying $B$-structures of the
trivial $G$-bundle on $\C{D}$. Hence it will suffice to show that $f$ is an isomorphism. 
Note that both $\C{B}'$ and $\C{B}$ are
equipped with a natural action of $G_{\C{D}}$, that $f$ is $G_{\C{D}}$-equivariant, and that
$\C{B}$ is a homogeneous space for the action of $G_{\C{D}}$.
Therefore it will suffice to check that the schematic preimage 
$\C{C}:=f^{-1}(o)$ of the point $o\in\C{B}$, corresponding to the standard $B$-structure,
consists of one reduced point.

By Cartan decomposition, $\C{C}_{\red}$ consists of one point $y_{\om}\in Gr(k)$. Explicitly, 
$y_{\om}=g_{\om}(y_0)$, where $y_0\in Gr(k)$ is a point corresponding to the trivial $G$-bundle on 
$\C{D}$, and $g_{\om}\in T(k((x)))\subset G(k((x)))=G_{\C{D}^*}(k)$ is the image of 
$x^{-1}\in k((x))\m=\B{G}_m(k((x)))$ under $\om:\B{G}_m\to T$. So it remains 
to show that the tangent space $T_{y_{\om}}(\C{C})$ is trivial.

For each dominant weight $\la$ of $G$, denote by $\C{L}_{0,\la}\subset triv_{\la}$ be the line subbundle
corresponding to the standard $B$-structure. Then $\C{C}$ is a schematic intersection inside $Gr$
of $Gr^0_{\om}$ with an ind-subscheme $\C{N}_{\om}$, consisting of points $(\C{G},\phi)$ such that 
$\C{L}_{0,\la}(\langle \la,\om\rangle x)$ is a line subbundle of $\C{G}_{\la}$ for each $\la$. 
Let $N$ be the unipotent radical of $B$. Then $\C{N}_{\om}$ is a homogeneous space for 
the action of the group 
$N_{\C{D}^*}$. Therefore we get a surjective map $p:N_{\C{D}^*}\to \C{N}_{\om}$ sending
$u$ to $u(y_{\om})$ which induces a surjection $dp:\Lie N(k((x)))\to T_{y_{\om}}(\C{N}_{\om})$.
Hence we have to check that  $dp(u)\notin T_{y_{\om}}(Gr_{\om})$ for each 
$u\in \Lie N(k((x)))\sm \Ker (dp)$. 

For each $u\in \Lie N(k((x)))$, denote by $\C{G}_u\subset G\times\C{D}^*\times\Spec k[t]/(t^2)$ the 
$G$-bundle on $\C{D}\times\Spec k[t]/(t^2)$,
corresponding to $dp(u)\in T_{y_{\om}}(Gr)\subset Gr(k[t]/(t^2))$. Then our assertion is 
equivalent to the fact that for each  $u\in \Lie N(k((x)))\sm \Ker (dp)$, there exists a 
dominant weight $\la$ of $G$ such that 
$(\C{G}_u)_{\la}\subset V_{\la}\times\C{D}^*\times\Spec k[t]/(t^2)$ is not contained in 
$(V_{\la}\times\C{D}\times\Spec k[t]/(t^2))(\langle \la,\om\rangle x)$.

Note that $\Lie N(k((x)))$ decomposes as the direct sum 
$\oplus_{\al}k((x))$, where $\al$ runs over the set $\Dt_+$ of all positive roots of $G$, and that
$\Ker (dp)=\oplus_{\al}x^{-\langle \om,\al\rangle }k[[x]]$.  For each 
$u\in\Lie N(k((x)))$ and $\al\in\Dt_+$, denote by $u_{\al}\in k((x))$ the $\al$-component of $u$.
Choose a basis $\{v_i\}_i$ of $V_{\la}$ consisting of $T$-eigenvectors and denote by 
$\mu_i\in X^*(T)$ the weight of $v_i$. Then 
$\{x^{-\langle \mu_i,\om\rangle }(v_i+t\sum_{\al}u_{\al}\al(v_i))\}_i$ generate 
the vector bundle $(\C{G}_u)_{\la}$.
Therefore we have to check that for each 
$u\in \Lie N(k((x)))\sm \Ker (dp)$ there exists $\la\in X_*^+(T)$, a $T$-eigen vector $v$ of $V_{\la}$ 
of weight $\mu$, and $\al\in\Dt_+$ such that $\al(v)\neq 0$ and 
$u_{\al}\notin x^{-\langle \la-\mu,\om\rangle }k[[x]]$. 
As $u\notin\Ker (dp)$, there exists $\al$ such that $u_{\al}\notin x^{-\langle \al,\om\rangle }k[[x]]$. 
Thus it remains to show that for each $\al\in\Dt_+$, there exists $\la\in X_*^+(T)$ and a 
$T$-eigenvector $v$ of $V_{\la}$ of weight $\la-\al$ such that $\al(v)\neq 0$. 
 
Using the representation theory of $SL_2$, the last statement is equivalent 
to the assertion that for each $\al$ there exists $\la$ such that 
$\langle\check{\al},\la\rangle$ (where $\check{\al}$ is the coroot corresponding to $\al$)
is not divisible by $\har k$. Assume that this is not the case.
Then for each $\mu\in X^+_*(T)$, the product $\langle\check{\al},\mu\rangle$ is divisible by
$\har k$. Since $\langle\check{\al},\al\rangle=2$, we thus get that $\har k=2$.
Let $G_1,\ldots,G_l$ be all the simple factors of $G^{\ad}$ numbered in a way that $\al$ is a 
root of $G_1$. Replacing $\al$ by a Weyl group conjugate, we can assume that $\al$ is simple.
As for every other simple root $\beta$ of $G_1$, we have $\langle\check{\al},\beta\rangle$ is even,
we see from the classification of simple groups that $G_1$ is either $PGL_2$ or $PO_{2m+1}$. 

It remains to show that $G_1$ is a direct factor of $G$. Since $G_1$ is adjoint, and the center
of $G_1^{\ssc}$ consists of two elements, it remains to show that the canonical homomorphism 
$G_1^{\ssc}\to G$ is not injective. To see this, denote by $T_1\subset G^{\ssc}_1$ the preimage of $T$. 
Then $X_*(T_1^{\ssc})$ is generated by coroots of $G_1^{\ssc}$, 
hence there exists $\mu\in X^*(T_1^{\ssc})$ such that $\langle\check{\al},\mu\rangle=1$. Therefore
our assumption on $\al$ implies that the restriction map $X^*(T)\to  X^*(T_1^{\ssc})$ is not surjective,
or, what is the same, the homomorphism $G_1^{\ssc}\to G$ is not injective, as claimed.

\begin{Rem} \label{R:reduced}
a) By the above argument, the condition in the proposition is not only sufficient 
but also necessary, that is, $Gr^0_{n,\ov{\om}}$ is not reduced if $k$ is a field of characteristic two,
and $G$ has a direct factor isomorphic to $PGL_2$ or $PO_{2m+1}$.

b) It would be interesting to check whether the full stack $Gr_{\om}$ is reduced.
This seems to be the case at least in some simple cases (e.g. for $G=SL_2$).
A positive answer to this question would imply that both $FBun_{n,\ov{\om}}$ 
and  $Hecke_{n,\ov{\om}}$ are reduced. In particular, a variant of $Hecke_{n,\ov{\om}}$, 
considered in \rr{inv} d), would then coincide with the original one.  
\end{Rem}

b) We start from showing that $Gr_{\om}$ is irreducible. As
$Gr^0_{\om}$ is a homogeneous space for the action of a connected
group scheme $G_{\C{D}}$, it is irreducible. Hence it will suffice to show that for
every two dominant coweights satisfying $\la_1<\la_2$, the orbit corresponding
to $\la_1$ lies in the closure of that of $\la_2$. For this we may
assume that the difference $\la_2-\la_1$ is a positive coroot
$\al$ of $G$. Indeed, by the lemma of Stembridge (see e.g., \cite[Lem
2.3]{Ra}), there exists a sequence of dominant coweights
$\la_1=\mu_0<\mu_1<\ldots<\mu_r=\la_2$ such that any two
neighboring $\mu$'s differ by a positive coroot. Next we may
assume that $G$ is of semi-simple rank one. Indeed, choose a maximal torus $T$
of $G$, and let $G'$ be the subgroup of $G$ generated by $T$
together with the image of the canonical morphism $SL_2\to G$,
corresponding to $\al$. Then $G'$ is a reductive group of semi-simple rank one, and the statement for $G'$ 
implies that for $G$. As any
reductive group of semi-simple rank one is a product of a torus with either
$GL_2$, $SL_2$ or $PGL_2$, it remains to check the statement for
$GL_2$. In this case the irreducibility of $Gr_{\om}$
easily follows from explicit resolution of singularities.

\begin{Rem} \label{R:irr}
The irreducibility of $Gr_{\om}$ is essentially equivalent (using Bott--Samelson
resolution of singularities) to the fact that the standard order on coweights is
induced by the Bruhat order on the affine Weyl group.
\end{Rem}

Next, as in the beginning of a), the assertion for $Gr'_{n,\ov{\om}}$ follows immediately from that 
for $Gr_{\om}$. In particular, we get
the statement for $Gr_{1,{\om}}$. The assertion for $Gr_{n,\ov{\om}}$ will be shown by induction on $n$. 
First by \rl{red} c) we get the irreducibility of
$Gr_{n,\overline{\om}}\times_{X^n} (X^n\sm\Dt)$. Denote now by $\ov{Gr}_{n,\overline{\om}}$ the closure of 
$Gr_{n,\overline{\om}}\times_{X^n} (X^n\sm\Dt)$ in  $Gr_{n,\overline{\om}}$.
For each $i\neq j$, consider the closed subscheme $(\ov{Gr}_{n,\overline{\om}})_{|x_{i}=x_j}$ of 
$\ov{Gr}_{n,\overline{\om}}$. As it is given locally by one equation in $\ov{Gr}_{n,\overline{\om}}$, 
it is of codimension one. On the other hand, by induction hypothesis, 
$(Gr_{n,\overline{\om}})_{|x_{i}=x_j}$ is irreducible and has the same dimension as 
$(\ov{Gr}_{n,\overline{\om}})_{|x_{i}=x_j}$. Thus $(\ov{Gr}_{n,\overline{\om}})_{|x_{i}=x_j}$ is 
contained in  $Gr_{n,\overline{\om}}$ for each $i\neq j$. Hence $\ov{Gr}_{n,\overline{\om}}=Gr_{n,\overline{\om}}$, as claimed.
\end{proof}

\begin{Lem} \label{L:small}
The forgetful morphism $\pi:Hecke'_{n,\overline{\om}}\to Hecke_{n,\overline{\om}}$ is
small.
\end{Lem}
\begin{proof}
As the statement is well known to experts, we will just sketch the argument for the 
convenience of the reader.
Observe first that it follows from basic properties of Coxeter groups that the Bott--Samelson resolution 
$\wt{Gr_{\om}}\to Gr_{\om}$ is semi-small. Our statement is a formal consequence
of this fact. Indeed, consider the stratification of $X^n$, given by diagonals $x_i=x_j$. As the restriction
of $\pi$ to the open stratum of $X^n$ is an isomorphism (by \rl{repr}), it will suffice to show 
that the restriction of $\pi$ to each stratum is semi-small. Moreover, by \rl{red} a) and the 
induction hypothesis, we have to check the assertion only for the closed stratum  $x_1=\ldots=x_n$.
Furthermore, by \rl{red} c),d), it will suffice to show that for each $(x,\C{G})\in X\times Bun$ the restriction $\pi_y$
of $\pi$ to $y:=(x,\ldots,x;\C{G})\in X^n\times Bun$ is semi-small. The last assertion follows from the observations 
that the fiber of $Hecke_{n,\overline{\om}}$ over $y$ is isomorphic to  $Gr_{\om_1+\ldots+\om_n}$, and  
 the (semi-small) Bott--Samelson resolution of $Gr_{\om_1+\ldots+\om_n}$ factors though $\pi_y$.
\end{proof}

\begin{Prop} \label{P:tate}
The restriction of the IC-sheaves of $Hecke_{n,\overline{\om}}$ and $Hecke'_{n,\overline{\om}}$
to each stratum isomorphic to a direct sum of
complexes of the form $\overline{\B{Q}_l}(k/2)[k]$ with the parity of
$k$ is the same as that of $\dim  Hecke_{n,\overline{\om}}=\dim  Hecke'_{n,\overline{\om}}$.
\end{Prop}
\begin{proof}
First we will show the corresponding statement for $Gr_{\om}$.
Let (the Iwahori subgroup) $I\subset G_{\C{D}}$ be the preimage of $B\subset G$ under the natural 
projection $G_{\C{D}}\to G$.
Since the IC-sheaf of  $Gr_{\om}$ is $G_{\C{D}}$-equivariant and since each stratum of $Gr_{\om}$
has an open $I$-orbit, it remains to
show the corresponding statement for the restriction of the IC-sheaf
to each $I$-orbit. For this we will use the same strategy as in
\cite[A.7]{Ga}. Consider the Bott--Samelson resolution
$\pi:\wt{Gr_{\om}}\to Gr_{\om}$. By the decomposition theorem,
$\IC_{Gr_{\om}}$ is a direct summand of
$\pi_!(\overline{\B{Q}_l})(\dim Gr_{\om})[2\dim Gr_{\om}]$. Therefore it will suffice to show
that the restriction of  $\pi_!(\overline{\B{Q}_l})$ to each
$I$-orbit  is a direct sum of complexes of the form
$\overline{\B{Q}_l}(k)[2k]$. Consider the stratification of
$\wt{Gr_{\om}}$ by $I$-orbits. As $I$ is pro-unipotent, each
stratum of  $\wt{Gr_{\om}}$ is an $\B{A}^N$-bundle (for some $N$)
over the corresponding stratum of $Gr_{\om}$. By the proper base
change theorem, it will therefore suffice to show the statement for
each fiber.

Thus we are reduced to showing that if $\rho:X\to\fq$ has a
stratification by  affine spaces, then $\rho_!(\overline{\B{Q}_l})$
decomposes as a direct sum of complexes of the form
$\overline{\B{Q}_l}(k)[2k]$. Let $N$ be the dimension of $X$, and let $U$ be the
disjoint union of (open) strata of $X$ of top dimension. By
induction, we may assume that the statement holds for $X\sm U$. As
the statement clearly holds for affine spaces, hence for $U$, it would suffice 
to show that $\rho_!(\overline{\B{Q}_l})$ decomposes as a direct
sum $\rho_!(\overline{\B{Q}_l}_{|U})\oplus\rho_!(\overline{\B{Q}_l}_{|X\sm U})$. 
As the statement is equivalent to splitting of the canonical distinguished triangle
$\rho_!(\overline{\B{Q}_l}_{|U})\to\rho_!(\overline{\B{Q}_l})
\to\rho_!(\overline{\B{Q}_l}_{|X\sm U})\to$, the assertion for $Gr_{\om}$ now follows from the 
fact that $\rho_!(\overline{\B{Q}_l}_{|U})=\tau_{\geq 2N}\rho_!(\overline{\B{Q}_l})$.

In the global case, we will show our assertion by induction on $n$. For $n=1$, it is an 
immediate consequence of the case of $Gr_{\om}$ was proved above (use \rl{red} c),d)). As
$Hecke'_{n,\overline{\om}}=Hecke_{1,\om_1}\times_{Bun} Hecke_{1,\om_2}\times_{Bun}\ldots\times_{Bun} Hecke_{1,\om_n}$,
the assertion for  $Hecke_{1,\om}$ implies that for $Hecke'_{n,\overline{\om}}$.
For $n>1$, take any stratum $S$ of
$Hecke_{n,\overline{\om}}$. We have two cases: either $S$ lies
over $X^n\sm \Dt$, or $S$ is a stratum of some 
$(Hecke_{n,\overline{\om}})_{|x_{i}=x_j}$. In the first case, the
statement follows from that for $Hecke'_{n,\overline{\om}}$ and \rl{repr}.

In the second case, the statement will follow from the induction
hypothesis if we show that the restriction of the IC-sheaf of 
$Hecke_{n,\overline{\om}}$ to $(Hecke_{n,\overline{\om}})_{|x_{i}=x_j}\subset Hecke_{n-1}$
is a direct sum of IC-sheaves of closed strata. 
By \rl{small}, $\pi_!$ maps the IC-sheaf of $Hecke'_{n,\overline{\om}}$ to that of $Hecke_{n,\overline{\om}}$.
By the proper base change theorem, we are therefore reduced to the case $n=2$, which in its turn 
reduces to the case of $Gr_{2,\overline{\om}}$.
In this case, the assertion over $\overline{\fq}$ is shown in 
the proof of \cite[Prop 1]{Ga}, and the decomposition over $\fq$ easily follows from the
fact that each fiber of $Gr'_{2,\overline{\om}}\to Gr_{2,\overline{\om}}$ has a stratification by affine spaces
(compare the proof of \rl{small}). 
\end{proof}

\end{document}